\input amstex.tex
\mag=\magstep1
\documentstyle{amsppt}
\vcorrection{-0.6cm}
\topmatter
\title 
     Analysis on the minimal representation of $O(p,q)$
\\
     --  III.  ultrahyperbolic equations on $\Bbb R^{p-1,q-1}$
\endtitle
\author 
  Toshiyuki KOBAYASHI
  and 
  Bent \O RSTED
\endauthor
\affil 
RIMS Kyoto and SDU-Odense University
\endaffil
\address
  RIMS, Kyoto University,
   Sakyo-ku, Kyoto, 606-8502, Japan
\newline\indent
  Department of Mathematics and Computer Science,
  SDU-Odense University, Campusvej 55, DK-5230, Odense M, Denmark
\endaddress
\email{toshi\@kurims.kyoto-u.ac.jp;
       orsted\@imada.sdu.dk}\endemail
\abstract
For the group $O(p,q)$ we give a new construction of its minimal
unitary representation via Euclidean Fourier analysis. This is
an extension of the $q = 2$ case, where the representation is the
mass zero, spin zero representation realized in a Hilbert space of
solutions to the wave equation.  
The group $O(p,q)$ acts as the M\"obius group of conformal transformations
 on $\Bbb R^{p-1, q-1}$,
 and preserves a space of solutions of 
 the ultrahyperbolic Laplace equation on $\Bbb R^{p-1, q-1}$.
We construct in an intrinsic and natural way a Hilbert space
 of solutions so that $O(p,q)$ becomes a continuous irreducible
 unitary representation in this Hilbert space.
We also prove that this representation is unitarily equivalent to
 the representation on $L^2(C)$, where $C$
is the conical subvariety of the nilradical of a maximal parabolic
subalgebra obtained by intersecting with the minimal nilpotent
orbit in the Lie algebra of $O(p,q)$. 

\endabstract
\endtopmatter
\NoRunningHeads

\mag=\magstep1
\def \sgn{\operatorname{sgn}}
\def\xabste{\cite{1}}
\def \xbail{\cite{2}}
\def \xbz{\cite{3}}
\def\Sahi{\cite{4}}
\def \xerdHigI{\cite{5}}
\def \xerdIntII{\cite{6}}
\def \xgs{\cite{7}}
\def \GSt{\cite{8}}
\def\xHo{\cite{9}}
\def\xkdecomp{{\cite{10}}}
\def\xkorsI{\cite{11}}
\def\xkorsII{\cite{12}}
\def \xkos{\cite{13}}
\def\xtodorov{\cite{14}}
\def\trF{{ {}_2 \! F_1}}
\def\spr#1#2{{\Cal H^{#1}(\Bbb R^{#2})}}
\def \trans{{}^t \!}

\def\diag{\operatorname{diag}}
\def \tilLap#1{{\widetilde {\Delta}}_{#1}}
\def \sec#1{{\vskip 0.5pc\noindent$\underline{{\bold{\ch.\sc.}} \; \text{{#1}}}$\enspace}}
\def\prince#1#2#3{{#1}\text{-}\operatorname{Ind}_{\Pmax}^G({#3} \otimes \Bbb C_{#2})}

\def \Image{\operatorname{Image}}
\def \Ker{\operatorname{Ker}}
\def \num{{\ch.\sc}} %
\define \set#1#2{\{{#1}:{#2}\}}
\def\Mmax{M^{\roman{max}}}
\def\Amax{A^{\roman{max}}}
\def\Nmax{N^{\roman{max}}}
\def\Pmax{P^{\roman{max}}}
\def \trans{{}^t \!}
\def\rarrowsim{\smash{\mathop{\,\rightarrow\,}\limits
  ^{\lower1.5pt\hbox{$\scriptstyle\sim$}}}}
\def\larrowsim{\smash{\mathop{\,\leftarrow\,}\limits
  ^{\lower1.5pt\hbox{$\scriptstyle\sim$}}}}
\def \Vpq{V^{p,q}}
\def \pmi#1{{#1}_{\pm}^{(i)}}
\def\tilPsi{{\Psi^*_{\frac{n-2}2}}} %
\def\invtilPsi{{(\Psi^*_{\frac{n-2}2, \epsilon})^{-1}}}
\def\minflat{{\varpi^{\roman{min}}_{\Bbb R^{p-1,q-1}}}} 
\def\ch{0}
\def\sc{1}
\define\zprime{z_1, \dots, \widehat{z_i}, \dots, z_n}
\define\zetaprime{\zeta_1, \dots, \widehat{\zeta_i}, \dots, \zeta_n}
\def \gk{$(\frak g, K)$}

\def\nbar#1{\overline{n}_{#1}}
\def\fLap{\square_{\Bbb R^{p-1,q-1}}} 
\document

\centerline{{\bf {Contents}}}
\smallskip
\item{ \S 1.}\enspace
 Introduction
\item{ \S 2.}\enspace
 Ultrahyperbolic equation on $\Bbb R^{p-1, q-1}$ and conformal group
\item{ \S 3.}\enspace
 Square integrable functions on the cone
 \item{ \S 4.}\enspace
 Green function and inner product
\item{ \S 5.}\enspace
   Bessel function and an integral formula of spherical functions 
\item{ \S 6.}\enspace
 Explicit inner product on solutions $\fLap f =0$

\def\ch{1}
\head
\S \ch.  Introduction
\endhead
\def\sc{1}
\sec{}
In this paper we study the symmetries of the ultrahyperbolic
Laplace operator on a real finite-dimensional vector space
equipped with a non-degenerate symmetric bilinear form. 
We shall work in coordinates so that the operator becomes 
$$
\fLap 
 \equiv   \square_z
 := \frac{\partial^2}{{\partial z_1}^2} + \dots +
    \frac{\partial^2}{{\partial z_{p-1}}^2}
    - \frac{\partial^2}{{\partial z_p}^2} - \dots -
    \frac{\partial^2}{{\partial z_{p+q-2}}^2},
$$  
on $\Bbb R^n = \Bbb R^{p-1, q-1}$.      
In the case of 
Minkowski space ($q = 2$) we are studying the wave equation, which
is well-known to have a conformally invariant space of 
solutions, see \xtodorov.  
This corresponds to the fact that the equation $\fLap f = 0$
in this case describes a particle of zero mass. 
Incidentally, it may also be interpreted as the bound states of the
Hydrogen atom, namely each energy level corresponds to a $K$-type -
for $(p,q) = (4, 2)$.
This gives the angular momentum values by further restriction to $O(3)$.
In general
the indefinite orthogonal group $G = O(p,q)$ acts as the M\"o{}bius
group of meromorphic conformal transformations on $\Bbb R^{p-1,q-1}$,
leaving a space of solutions to $\fLap f = 0$ invariant.

\def\sc{2}
\sec{}
The main object of the present paper is to construct in an intrinsic
and natural way a Hilbert space of solutions of $\fLap$ so that the action
of $O(p,q)$ becomes a continuous {\bf unitary} irreducible representation 
in this Hilbert space for $(p,q)$ such that $p, q \ge 2$ and $p+q>4$ is even.
{}From an algebraic view point of representation theory,
 our representations are:
\newline\indent
i) minimal representations if $p + q \ge 8$ (i.e. the annihilator is the Joseph ideal).
\newline\indent 
ii) {\bf not} spherical if $p\neq q$ (i.e. no non-zero $K$-fixed vector).
\newline\indent
iii) {\bf not} highest weight modules of $SO_0(p,q)$ if $p, q \ge 3$.

In a long history of representation theory of semisimple Lie groups,
 it is only quite recent that our representations for $p, q \ge 3$
  have been paid attention,
  especially as minimal unitary representations;
they were first discovered by Kostant \xkos\ for $(p,q) = (4,4)$
 and generalized by Binegar-Zierau \xbz\ 
 as subrepresentations of degenerate principal series representations.
There is also another algebraic approach to the same representations
 by using the theta correspondence for the trivial representation 
of $SL(2, \Bbb R)$ by Huang-Zhu. Our previous papers 
\xkorsI{} and \xkorsII{} treated  the same representation
by geometric methods and with other points of view.
We think that such various approaches reflect a rich structure of 
 the minimal representations.

It is perhaps of independent interest that the (in some sense maximal
 group of) symmetries $O(p,q)$ of the space of solutions of $\fLap f = 0$
 lead to such a natural Hilbert space.
Our inner product $(\ ,\ )_W$ defined by an integration 
 over a non-characteristic hyperplane (see (\ch.5.1))
 is a generalization of the one coming
from energy considerations in the case of wave equations, and
even the translation invariance of the inner product contains some
new information about solutions.

It is also of independent interest from the representation theory 
 of semisimple Lie groups that our representations are
 unitarily equivalent to the representations on $L^2(C)$,
  where $C$ is the null cone of the quadratic form on $\Bbb R^{p-1, q-1}$.
This result is proved  via the Fourier transform in Theorem~4.9.
Such $L^2$-realizations of \lq\lq unipotent representations"
 is expected from the philosophy of the Kostant-Kirillov orbit method,
  but has not been proved except for some special cases of
  highest weight modules or spherical representations.
  
 We have avoided most of the
references to the theory of semisimple Lie groups and representation 
theory, and instead given direct constructions of the key objects,
such as for example the minimal $K$-type; this is given as an
explicit hypergeometric function, and we also calculate its Fourier
transform in terms of a Bessel function. 
 By application of explicit
differential operators forming the Lie algebra of $G$ we can generate
the whole Hilbert space of solutions beginning from the minimal $K$-type.

\def\sc{3}
\sec{}
For $q = 2$ (or $p = 2$) we are dealing with highest weight representations
(when restricting to the identity component $SO_0(p,2)$), 
 and these have been studied by many authors,
in particular in the physics literature. For a nice
introduction to this representation and its construction via geometric 
 quantization
(and more) see \GSt. 
In this case the $K$-types may be identified
with energy levels of the bound states of the Hydrogen atom, 
 and the smallest one with the bound state of lowest energy. 

We can summarize 
 the situation, covering both the classical Kepler
problem and its quantization in case $q=2$, as in the diagram below. Here
the left-hand side represents the classical descriptions of
respectively the Kepler problem and geodesic flow on the sphere;
by \lq\lq symplectic transform" we are alluding to the change of variables
between these two Hamiltonian systems as presented in \GSt. The
right-hand side involves the quantizations of these two systems, where
the wave-equation is considered
as the quantization of geodesic flow, also to be thought of as geometric optics. 
The quantum analogue of the \lq\lq symplectic transform" involves the Fourier transform. 
Finally we invoke conformal geometry and combine it with the Fourier transform,    
which in a different (and new as far as constructing Hilbert
spaces and unitary actions) way appears in passing from
the wave equation to the Fourier realization of solutions - this
is the last arrow on the right-hand side.    
$$
  \align
\text{Classical} &\Rightarrow \text{Quantum}
\\
{}
\\
{{
 \vbox{
         \offinterlineskip
        \hrule
            \halign
                 {&\vrule #   & \strut \quad \hfil # \; \cr
                  & Kepler problem    &                        \cr
                  }
            \hrule}}}^{\;\underset{\tsize{\Rightarrow}}\to {}}
 &{\vbox{
         \offinterlineskip
         \hrule
         \halign
                 {&\vrule # & \strut \quad \hfil # \; \cr
                  & Hydrogen atom \hphantom{}&\cr
}
\hrule}}
\\
\Downarrow
\text{symplectic transform}
&
\qquad\qquad\Downarrow\text{Fourier transform}
\\
{\vbox{
\offinterlineskip
\hrule
\halign
{&\vrule # & \strut \quad \hfil # \; \cr
&geometric optics &\cr
\cr}
\hrule}}
^{\;\underset{\tsize{\Rightarrow}}\to {}}
&{\vbox{
\offinterlineskip
\hrule
\halign
{&\vrule # & \strut \quad \hfil # \; \cr
&(A) wave equation &\cr
}
\hrule}}
\\
&\qquad\qquad\quad\text{conformal geometry}
\\
&\qquad\qquad\Downarrow\text{Fourier transform}
\\
&{\vbox{
\offinterlineskip
\hrule
\halign
{&\vrule # & \strut \quad \hfil # \; \cr
& (B) other realizations and \hphantom{m} &\cr
& explicit inner products &\cr
 & of minimal representations &\cr
}
\hrule}}
\endalign
$$       
The main focus of this paper is on the boxes (A) and (B). 
In particular,
 we give an explicit inner product in the model (A) 
 (Theorems~\ch.4 and \ch.5)
 and
 construct via Fourier transform a new realization of the minimal
 representation (Theorem~\ch.6)
 for general $p, q$.

\def\sc{4}
\sec{}
{}From now, suppose that $n := p+ q -2$ is an even integer greater
 than $2$, and $p, q \ge 2$.
Let us briefly state some of our main results in a more explicit way.

First, we find a formula of the Green function $E_0$
 of the ultrahyperbolic Laplace operator  $\fLap$,
 in Proposition~4.2,
 namely, $E_0$ is given by a constant multiple of the imaginary part of
 the regularized Schwartz distribution:
$$
  e^{\frac{\sqrt{-1}\pi(q-1)}2}
  ({x_1}^2+\dots + {x_{p-1}}^2-{x_p}^2-\dots-{x_{p+q-1}}^2
  +\sqrt{-1}0)^{1-\frac{n}2}.
$$
See the recent paper of H\"ormander \xHo{} for further details on
distributions associated with this ultrahyperbolic equation.
Then we construct solutions of $\fLap f = 0$
 by the integral transformation:
$$
  S : C^\infty_0(\Bbb R^n) \to C^\infty(\Bbb R^n),
 \quad   \varphi \mapsto E_0 * \varphi
\qquad
\text{(see (4.3.1)).}
$$ 
The image $S(C^\infty_0(\Bbb R^n))$ turns out to be \lq\lq large\rq\rq\
 in $\Ker \fLap$ (see \S 4.7, Remark~(2)).
On this image, we define
 a Hermitian form $(\ , \ )_N$ by
$$
  (f_1, f_2)_N := 
 \int_{\Bbb R^n} \int_{\Bbb R^n} E_0(y-x) \varphi_1(x)
 \overline{\varphi_2(y)} \ d x  d y,
\tag \num.1
$$
 where $f_i = E_0*\varphi_i$ $(i=1,2)$.
Here is a part of Theorem~4.7, which is the first of our main results:
\proclaim{Theorem~\num}
 $(\ , \ )_N$ is positive-definite on the image of $S$.
Furthermore,
 $O(p,q)$ acts as an irreducible unitary representation
 on its Hilbert completion $\Cal H$.
\endproclaim
We shall write $(\minflat, \Cal H)$ for this unitary representation.
We also prove that
 this representation is isomorphic to the minimal representation of $O(p,q)$,
 constructed previously by Kostant, Binegar-Zierau (\xkos, \xbz).
Thus, Theorem~1.4 may be regarded as a realization of the minimal 
 representation (with an explicit inner product)
 in the solution space of the ultrahyperbolic equation.

\def\sc{5}
\sec{}
The above definition of the inner product $(\ , \ )_N$ (see (\ch.4.1))
 uses the integral expression
 of  solutions of $\fLap f = 0$.
Can we write the inner product without knowing the preimage?
Yes, the second of our main results is to
 give an intrinsic inner product on the same solution space
  by using the Cauchy data.
For simplicity, we take $z_1 = 0$ as a non-characteristic
 hyperplane.
Then, we decompose a solution 
$$
f  = f_+ + f_-
$$
 such that $f_\pm(z_1, \dots, z_n)$ is holomorphic with respect to the first variable
 $z_1$ in the complex domain of 
 $\set{z_1 \in \Bbb C}{\pm \operatorname{Im}z_1>0}$
 of the $z_1$-variable.
This is an expression of $f$ as a hyperfunction,
 and such a pair $(f_+, f_-)$ can be obtained
  by the convolution in the $z_1$-variable (see (6.2.3)):
 $$
  f_\pm (z) = \frac{1}{2 \pi \sqrt{-1}} \cdot
                \frac{\pm 1}{z_1 \mp \sqrt{-1} 0}
                     * f(z_1, \dots, z_n),
$$
where the integration makes sense for $f$ with suitable decay at infinity.
Then we define a Hermitian form
$$
  (f, f)_W := \frac{1}{\sqrt{-1}} \int_{\Bbb R^{n-1}}
   \left(f_+\overline{\frac{\partial f_+}{\partial z_1}}
   -f_-\overline{\frac{\partial f_-}{\partial z_1}}\right)|_{z_1=0}
  \ d z_2 \dotsb d z_n.
\tag \num.1
$$
Then we shall prove that $(\ , \ )_W$ is independent of
 the specific choice of a non-characteristic hyperplane,
 as follows from the (non-trivial) isometric invariance.
Much more strongly, $(\ , \ )_W$ is conformally invariant.
A precise formulation for this is given in Theorem~6.2,
 which includes:
\proclaim{Theorem~\num}
$ 4 \pi (\ ,\ )_W = (\ ,\ )_N$.
In particular,
 $( \ , \ )_W$  is positive definite and $O(p,q)$-invariant.
\endproclaim
Hence,
 in place of Theorem~\ch.4,
 we can obtain the same irreducible unitary representation of $O(p,q)$
 on the Hilbert completion of a space of solutions 
 with respect to the inner product $(\ , \ )_W$.

An interesting property of this inner product is its large invariance group.
Even in the case of the usual wave equation ($q=2$ case)
 our approach to the Hilbert space of solutions and the 
 corresponding representation offers some new points of view.
In this case, %
 if we take 
 the non-characteristic hyperplane as fixed time coordinate 
 (namely, if $z_1$ stands for the time), 
 then the translational invariance amounts to a remarkable
 \lq\lq conservation law\rq\rq.
Instead, we can take  
 the non-characteristic hyperplane by fixing one of the space coordinates,
 and an analogous integration over the hypersurface (containing
 the time coordinate) still gives the same inner product !
As a final remark in \S 6.7, we note the connection to
the theory of conserved quantities for the wave equation ($q = 2$ case), such
as the energy and others obtained by the action of the conformal group.  

\def\sc{6}
\sec{}
The Gelfand-Kirillov dimension of our representation
 $(\minflat, \Cal H)$ is $p+q-3$.
So, we may expect that the representation could be
 realized on a $(p+q-3)$-dimensional manifold.
For this purpose,
 we define the null cone of the metric as
$$
   C := \set{\zeta \in \Bbb R^n}{{\zeta_1}^2+\dots+{\zeta_{p-1}}^2
         -{\zeta_p}^2-\dots-{\zeta_n}^2 =0}.
$$

The third of our main results is another realization of 
 the unitary representation $(\minflat, \Cal H)$ in
 a function space on a $(p+q-3)$-dimensional manifold $C$.
The Fourier transform $\Cal F$
 maps solutions of $\fLap f = 0$ to distributions supported on the null cone $C$.
Surprisingly, the inner product of our Hilbert space turns out to be
 simply the $L^2$-norm on $C$
 with respect to a canonical measure $d \mu$ (see (3.3.3)) !
Here is a part of Theorem~4.9:
We regard $L^2(C)$ as a subspace of distributions
 by a natural injective map $T : L^2(C) \to \Cal S'(\Bbb R^n)$.
\proclaim{Theorem~\num}
 $(2 \pi)^{-\frac{n}2} T^{-1}\circ \Cal F$
 is a surjective unitary operator from $\Cal H$ to $L^2(C)$.
\endproclaim

Theorem~\num\ defines an irreducible unitary representation of $G = O(p,q)$
 on $L^2(C)$, denoted by $\pi$, 
 which is unitarily equivalent to $(\minflat, \Cal H)$.
Since the maximal parabolic subgroup
 $\overline{\Pmax}$ of $G$ (see \S 2.7) acts on $\Bbb R^{p-1, q-1}$
 as affine transformations,
 the restriction $\pi|_{\overline{\Pmax}}$ has a very simple form,
 namely,
 the one obtained by the classical Mackey theory (see (3.3.5)).
In this sense, 
 Theorem~\num\ may be also regarded as an {\bf extension theorem}
 of an irreducible unitary representation
 from the maximal parabolic subgroup
 $\overline{\Pmax}$ to the whole group $G$.

\def\sc{7}
\sec{}
The fourth of our main results is about
 the representation $(\pi, L^2(C))$ as \gk-modules
 on the Fourier transform side,
 especially to find an explicit vector in the minimal $K$-type.

In the realization on $L^2(C)$,
 the action $\pi(g)$ is not simple to describe
 except for $g \in\overline{\Pmax}$.
Instead, we consider the differential action $d \pi$ of 
 the Lie algebra $\frak g_0$ on smooth vectors of $L^2(C)$,
 which turns out to be
 given by differential operators at most of second order (see \S 3.2).
This makes the analogy with the metaplectic representation
(where $G$ is replaced by the symplectic group) a good one
. Here we are
recalling the fact, that the even part of the metaplectic
representation may be realized as an $L^2$-space of functions
on the cone generated by rank one projections in $\Bbb R^n$.   

Moreover,
 by using a reduction formula of an Appell hypergeometric function,
 we find explicitly the Fourier transform of a Jacobi function
 multiplied by some conformal factor
 which equals to a scalar multiple of
$$
\psi_{0,e}(\zeta) :=
|\zeta|^{\frac{3-q}2} K_{\frac{q-3}2}(2 |\zeta|) d \mu \in \Cal S'(\Bbb R^n).
$$
Here
 $K_\nu(\zeta)$ is a modified Bessel function of the second kind.
This vector
 $\psi_{0,e}(\zeta)$
 corresponds to the bound state of lowest energy for $q=2$ case
. 
For general $p, q$,
 the $K$-span of 
 $\psi_{0,e}(\zeta)$
 generates the minimal $K$-type in the realization on $L^2(C)$.

We define a subspace $U$ of $\Cal S'(\Bbb R^n)$ to be the linear span 
 of its iterative differentials
$$
 d \pi (X_1) \dotsb d \pi(X_k) \psi_{0,e}(\zeta)
\qquad
 (X_1, \dots, X_k \in \frak g_0 \otimes_{\Bbb R} \Bbb C).
$$
What comes out of \S 5 may be formulated in this way
 (combining with Theorem 4.9, see \S 3.2 for notation): Suppose $p+q \in 2 \Bbb Z$, $p+q > 4$
 and $p \ge q \ge 2$.
\proclaim{Theorem~\num}\enspace
{\rm 1)}\enspace
 $|\zeta|^{\frac{3-q}2} K_{\frac{q-3}2}(2 |\zeta|)$ is a $K$-finite vector
 in $L^2(C)$.
\newline
{\rm 2)}\enspace
$U$ is an infinitesimally unitary
 $(\frak g, K)$-module via $\widehat\varpi_{\frac{n-2}2, \epsilon}$.
\newline
{\rm 3)}\enspace
 $U$ is dense in the Hilbert space $T(L^2(C))$.
\newline
{\rm 4)}\enspace
The completion of (2) defines an irreducible unitary representation
 of $G$ on $T(L^2(C))$, and then also on $L^2(C)$.
\endproclaim

In the paper \Sahi\ one finds a similar construction of Hilbert spaces
and unitary representations for Koecher-Tits groups associated with
semisimple Jordan algebras under the assumption that the representations
are spherical, and there also occur Bessel functions
as spherical vectors. In our situation the representations are {\bf not}
spherical if $p\neq q$.
Our approach is completely different from \Sahi,
 and
 contrary to what is stated in \Sahi\ (p\. 206)
 we show that for $G = O(p,q)$ it is possible to extend the Mackey
 representation of the maximal parabolic subgroup to the whole group.
Furthermore, even for $p=q$ case,
 our approach for Theorem~1.7 has an advantage that
 we give the exact constants normalizing the unitary correspondence
 between the minimal $K$-type 
 in other realizations and the Bessel function in our realization on $L^2(C)$
 (see Theorem~5.5).

\def\sc{8}
\sec{}
The paper is organized as follows:
We begin by recalling some results from conformal geometry and facts
about the conformal group, in particular
 in \xkorsI{}.
 In section 3 we give the basic setup for a realization 
 on the null cone via Fourier transform.
Then we construct the intertwining operator
from the minimal representation to the model treated here and calculate
the new expression for the inner product (see Theorem~1.4). 
We show in Proposition 4.2 that the
Green function of $\fLap$ has a Fourier transform equal to the
invariant measure on the null-cone, allowing one more expression for
the inner product (see Theorem~1.6); also we obtain from this an intertwining operator
from test functions to solutions.
Indeed, in section 4 Proposition 4.6 we prove that the Green function is up to
a constant exactly the kernel in the Knapp-Stein intertwining
integral operator between degenerate principal series representations
at the parameters we study; this makes the proof of unitarity of
the minimal representation elementary and 
formulated independent
of the semi-simple theory
usually employed here. Note that all
normalizing constants are computed explicitly.
Lemma 2.6 states the irreducibility and unitarizability,
which we use; we give in \xkorsII\ independent proofs of these facts.    

In section 5 we construct the lowest $K$-type as a modified Bessel function,
whose concrete properties are important for $K$-type information about
$L^2(C)$. 
The idea here is to use a classical formula on the Hankel transform due to Baily in 1930s,
 and then apply reduction formulae of an Appell hypergeometric function of two variables.

Section 6 contains formulae for the inner product $( \ , \ )_W$ in terms of
integration over a Cauchy hypersurface.  
Summarizing, we give five different realizations of the
inner product together with the normalizations of these
relative to each other. Namely, in addition to 
$( \ , \ )_N$ and $( \ , \ )_W$ we also define three more:
$( \ , \ )_M$ (coming from a pseudodifferential operator on $ M = S^{p-1} \times S^{q-1}$),
$( \ , \ )_A$ (coming from a normalized Knapp-Stein intertwining operator), and
finally $( \ , \ )_C$, which is just $L^2(C)$.    
This is seen in the key diagram (see Section 4.11)
$$
\alignat1
 C_0^\infty(\Bbb R^n)
 \overset{S}\to\rightarrow
 &\tilPsi(\tilLap{M})
 \overset{\Cal F}\to\rightarrow
 \Cal S'(\Bbb R^n)
 \overset{T}\to\hookleftarrow
 L^2(C)
\\
& \quad \cap
\\ &\Ker \fLap
\endalignat
$$ 
where the spaces correspond to four different ways of generating 
solutions to our ultrahyperbolic equation. $S$ will be an integral transform
 against the Green kernel (essentially, a
 Knapp-Stein intertwining operator with a specific parameter),
  and $\Cal F$ the Fourier
transform, mapping solutions to distributions supported on the null cone $C$. 
Correspondingly to the various ways of generating solutions, we write down explicitly the
unitary inner product and its Hilbert space.  
We have tried to avoid the use of any semi-simple
theory and stay within classical analysis on spheres and
Euclidean spaces; still our treatment may also be of
interest to people working with the classification of the
unitary dual of semi-simple Lie groups, since we are
providing new models of some unipotent representations.
Tools like the standard Knapp-Stein intertwining operators
become very natural to use here, also from the more elementary viewpoint,
and the close connection between these and Green functions
for ultrahyperbolic differential operators seems not to have
been noticed before. 

The first author expresses his sincere gratitude to SDU - Odense University 
for the warm hospitality.

\def\ch{2}
\head
\S \ch.  Ultrahyperbolic equation on $\Bbb R^{p-1, q-1}$ and
 conformal group %
\endhead
\def\sc{1}
\sec{}
As explained in the Introduction,
 we shall give a flat picture, the so-called
$N$-picture, of the minimal representation,
 which is connected to classical facts about conformal geometry in $\Bbb R^n$.
We shall give a unitary inner product in this realization (see Theorem 6.2)
 and also in its Fourier transform (Theorem 4.9),
 together with an explicit form of minimal $K$-type in this realization
 (see Theorem 5.5).

We shall assume $p+q \in 2 \Bbb N$,
 $p \ge 2$, $q \ge 2$ and $(p,q) \neq (2,2)$.
The parity condition $p + q \in 2 \Bbb N$ is not necessary
 when we consider
 a representation of the parabolic 
 subgroup $\overline{\Pmax}$ or of the Lie algebra $\frak g$.
Indeed, it will be interesting to relax this parity condition in order to
obtain an infinitesimally unitary representation, which does not
integrate to a global unitary representation of $G$.
 
Throughout this paper, we let 
$$
                   n = p+q-2
$$

This section is written in an elementary way,
 intended also for non-specialists of semisimple Lie groups.
\S \ch.2 and \S 2.6 review the needed results in \xkorsI.

\def\sc{2}
\sec{}
We recall some basic fact of the distinguished
 representation of a conformal group (see \xkorsI, \S~2).
Let $M$ be an $n$-dimensional manifold
 with pseudo-Riemannian metric $g_M$.  
We denote by $\Delta_M$
 the Laplace  operator on $M$,
 and by $K_M$ the scalar curvature of $M$.
The Yamabe operator is defined to be
$$
     \tilLap{M} := \Delta_M - \frac {n-2}{4(n-1)} K_M.
$$
Suppose $(M, g_M)$ and $(N, g_N)$ are pseudo-Riemannian manifolds.
A local diffeomorphism
 $\Phi \: M \to N$ is called
 a {\it{conformal map}}
 if there exists a positive-valued function $\Omega$ on $M$
 such that
$
     \Phi ^* g_N = \Omega ^2 g_M.  
$
For $\lambda \in \Bbb C$,
 we introduce a twisted pull-back
$$
   \Phi_\lambda^* \: C^\infty(N) \to C^\infty(M), 
                  \ f \mapsto \Omega^\lambda \cdot f\circ \Phi.
\tag \num.1
$$
Then the conformal quasi-invariance of the Yamabe operator is
expressed by:
$$
  \Phi_{\frac {n+2} 2}^* \tilLap{N} 
  = \tilLap{M} \Phi_{\frac {n-2}2}^*.  
\tag \num.2
$$
Let $G$ be a Lie group acting conformally on $M$.  
If we write the action as $x \mapsto L_h  x$ $(h \in G, x \in M)$,
 we have a positive function $\Omega(h,x) \in C^\infty(G \times M)$
 such that
$$
     L_h^* g_M = \Omega(h,\cdot)^2 \ g_M
     \qquad (h \in G).  
$$
We form a representation $\varpi_{\lambda}$
 of $G$, with parameter $\lambda \in \Bbb C$,
 on $C^{\infty} (M)$ as follows:
$$
     \varpi_{\lambda}(h^{-1}) f (x)= \Omega(h,x)^{\lambda} f (L_h x),
     \quad
     (h \in G, f \in C^{\infty}(M), x \in M).  
\tag \num.3
$$
 Note that the right-hand side is given by the twisted
pull-back $(L_h)_{\lambda}^*$ according to the notation (\num.1).
Then, Formula (\num.2) implies that
$\tilLap{M} \: C^{\infty} (M) \to C^{\infty} (M)$
 is a $G$-intertwining operator
 from $\varpi_{\frac {n-2} 2}$ to $\varpi_{\frac {n+2}2}$.  
Thus,
 we have constructed a distinguished representation of the conformal group:
\proclaim{Lemma~\num\ {\rm (see \xkorsI, Theorem~2.5)}}
 $\Ker \tilLap{M}$ is a representation space of
 the conformal group  $G$
 of a pseudo-Riemannian manifold $(M,g_M)$,
 through $\varpi_{\frac {n-2}2}$.  
\endproclaim
If $(N, g_N)$ is also a pseudo-Riemannian manifold on which the same group
 $G$ acts conformally,
 then one can also define a representation $\varpi_{\lambda, N}$
 on $C^\infty(N)$.
Then the twisted pull-back $\Phi_\lambda^*$ is a $G$-intertwining operator.

\def \sc{3}
\sec{}
Here is a setup on which we construct the minimal representation of
 $O(p,q)$ by applying Lemma~\ch.2.
Let $p, q \ge 2$.
We note
$
   n = p+q-2.
$
We write $\{e_0, \dots, e_{n+1}\}$ for a standard basis of $\Bbb R^{p+q}$ 
 and the corresponding coordinate as
$$
 (v_0, \dots, v_{n+1}) = 
  (x,y) = (v_0, z', z'', v_{n+1}),
$$
  where
 $x \in \Bbb R^p, y \in \Bbb R^q, z' \in \Bbb R^{p-1}, z'' \in \Bbb R^{q-1}$.
The notation $(x,y)$ will be used for $S^{p-1} \times S^{q-1}$, 
 while $(z', z'')$ for $\Bbb R^n = \Bbb R^{(p-1)+(q-1)}$.
The standard norm on $\Bbb R^l$
 will be written as $|\cdot|$ ($l = p-1, p, q-1, q$).

We denote by $\Bbb R^{p,q}$ the pseudo-Riemannian manifold $\Bbb R^{p+q}$
 equipped with the flat pseudo-Riemannian metric:
$$
 g_{\Bbb R^{p,q}} = {d v_0}^2 + \cdots + d {v_{p-1}}^2 - {d v_{p}}^2 
         - \cdots - d {v_{n+1}}^2.  
\tag \num.1
$$
We put two functions on $\Bbb R^{p+q}$ by
$$
\alignat2
           &\nu \: \Bbb R^{p+q} \to \Bbb R, &&\quad (x,y) \mapsto |x|,  
\tag \num.2
\\
           &  \mu \: \Bbb R^{p+q} \to \Bbb R, 
         \quad
         &&(v_0, \dots, v_{n+1}) \mapsto \frac12(v_0 + v_{n+1}).
\tag \num.3
\endalignat
$$
 and define three submanifolds of $\Bbb R^{p,q}$ by
$$
\alignat3
     &\Xi  && := \set{(x, y) \in \Bbb R^{p,q}}{|x|=|y|\neq 0}, &&
\\  
      &M  &&:= \set{v \in \Bbb R^{p,q}}{\nu(v) = 1} \cap \Xi
            &&   = S^{p-1} \times S^{q-1},  
\\
     &N &&:= \set{v \in \Bbb R^{p,q}}{\mu(v) = 1} \cap \Xi
        &&   \underset{\iota}\to{\larrowsim} \Bbb R^{n}. 
\endalignat
$$
where the bijection $\iota \: \Bbb R^n \to N$ is given by
$$
  \iota \: \Bbb R^{n} \to N,
  (z', z'') \mapsto
(1-\frac{|z'|^2 - |z''|^2}4, z', z'',
1+\frac{|z'|^2 - |z''|^2}4).
\tag \num.4
$$

We say a hypersurface $L$ of $\Xi$ is {\it transversal to rays}
 if the projection
$$
  \Phi\: \Xi \to M, 
   v \mapsto \frac{v}{\nu(v)}
\tag \num.5
$$
  induces a local diffeomorphism
 $\Phi|_L \: L \to M$.
Then, 
 one can define a pseudo-Riemannian metric $g_L$ 
 of signature $(p-1, q-1)$ on $L$ by the restriction of $g_{\Bbb R^{p,q}}$.
In particular,
 $M$ itself is transversal to rays,
 and the induced metric $g_{S^{p-1} \times S^{q-1}}$ equals
  $g_{S^{p-1}} \oplus (- g_{S^{q-1}})$,
 where $g_{S^{n-1}}$ denotes the standard Riemannian metric on the unit
 sphere $S^{n-1}$.
Likewise,
 the induced pseudo-Riemannian metric on $\Bbb R^n$ through
 $\iota \: \Bbb R^n \hookrightarrow \Bbb R^{p,q}$
 coincides with the standard flat pseudo-Riemannian metric
  $g_{\Bbb R^{p-1,q-1}}$ on $\Bbb R^n$.

\def\sc{4}
\sec{}
Let $I_{p,q} := \diag (1, \dots, 1, -1, \dots, -1) \in GL(p+q, \Bbb R)$.
The indefinite orthogonal group 
$$
  G = O(p,q) :=\set{g \in GL(p+q, \Bbb R)}{\trans{g} I_{p,q} g = I_{p,q}},
$$ 
 acts isometrically on $\Bbb R^{p,q}$
 by the natural representation,
 denoted by $v \mapsto g\cdot v$. %
This action stabilizes the light cone $\Xi$.
We note that
 the multiplicative group 
 $\Bbb R_+^{\times} := \set{r \in \Bbb R}{r >0}$ also acts on $\Xi$
 as dilation,
 which commutes with the linear action of $G$.
Then, using dilation,
 one can define an action of $G$ on $M $,
 and also a meromorphic action on $\Bbb R^{p-1, q-1}$
 as follows:
$$
\alignat3
 &L_{h} \: M\to M,
 \quad
 && v \mapsto \frac{h \cdot v}{\nu (h \cdot v)}
 &&\quad (h \in G),
\tag \num.1
\\
  &L_h \: \Bbb R^{p-1, q-1} \to \Bbb R^{p-1, q-1},
 \quad 
  &&z \mapsto
  \iota^{-1} \left( \frac{h \cdot \iota(z)}{\mu (h \cdot \iota(z))}\right)
 &&\quad
 (h \in G).
\tag \num.2
\endalignat
$$
Then,
 both of these actions are conformal:
$$
\alignat1
   (L_h)^* g_{M}
   &= \nu(h\cdot v)^{-2}  g_{M}
\tag \num.3
\\
   (L_h)^* g_{\Bbb R^{p-1,q-1}}
   &= \mu(h \cdot\iota(z))^{-2} g_{\Bbb R^{p-1,q-1}}.
\tag \num.4
\endalignat
$$
We note that (\num.2) and (\num.4) are
 well-defined if $\mu(h \cdot \iota(z)) \neq 0$.
In fact, $G$  acts only meromorphically on $\Bbb R^{p-1, q-1}$. 
An illustrative example for this feature
 is the linear fractional transformation of $SL(2, \Bbb C)$ on 
 $\Bbb P^1 \Bbb C = \Bbb C \cup \{\infty\}$,
 which is a meromorphic action on $\Bbb C$.
This example essentially coincides with (\num.2) for $(p,q) = (3,1)$,
 since $SL(2, \Bbb C)$ is locally isomorphic to $O(3,1)$
 and $\Bbb C \simeq \Bbb R^2$.

\def\sc{5}
\sec{}
The (meromorphic) conformal groups for the submanifolds $M$ and $N$ of $\Xi$
 are the same, namely, $G=O(p,q)$,
 while their isometry groups are different subgroups of $G$,
 as we shall see in Observation~2.5.
In order to describe them,
 we define subgroups $K$, $\Mmax$, $\Nmax$, $\Amax$ and $\overline{\Nmax}$
 of $G$ as follows:

First, we set
$$
\alignat2
 m_0 &:= -I_{p+q}, &&
\\
    K &:= G \cap O(p+q) &=  \ & \ O(p) \times O(q),
\\
    \Mmax_+ &:= \set{g \in G}{g \cdot e_0 = e_0, \ g \cdot e_{n+1} = e_{n+1}}
     \ & \ \simeq \ &O(p-1,q-1),
\\
    \Mmax &:= \Mmax_+ \cup m_0 \Mmax_+
    \ & \ \simeq &O(p-1,q-1) \times \Bbb Z_2.
\endalignat
$$
The Lie algebra of $G$ is denoted by $\frak g_0 = \frak{o}(p,q)$,
 which is given by matrices:
$$
   \frak g_0 \simeq \set{X \in M(p+q, \Bbb R)}{X I_{p,q} + I_{p,q} \trans{X} = O}.
$$
Next, we keep $n = p+q-2$ in mind and put
$$
\varepsilon_j = \cases
                       1 & (1 \le j \le p-1), 
                 \\
                      -1 & (p \le j \le n),
                \endcases
\tag \num.1
$$
and define elements of $\frak g_0$ as follows:
$$
\alignat3
   \overline{N}_j &:= E_{j, 0} + E_{j, n+1} - \varepsilon_j E_{0, j} 
            + \varepsilon_{j} E_{n+1, j} &&
\quad (1 \le j \le n),
\tag \num.2)(a
\\
   {N}_j &:= E_{j, 0} - E_{j, n+1} - \varepsilon_j E_{0, j} 
           - \varepsilon_{j} E_{n+1, j} &&
\quad (1 \le j \le n),
\tag \num.2)(b
\\
E &:= E_{0, n+1} + E_{n+1, 0}, &&%
\tag \num.2)(c
\endalignat
$$
 where $E_{i j}$ denotes the matrix unit.
Now, we define abelian subgroups of $G$ by
$$
    \overline{\Nmax}
          := \exp (\sum_{j=1}^n \Bbb R \overline{N}_j), 
   \quad
    \Nmax := \exp (\sum_{j=1}^n \Bbb R {N}_j),
   \quad
    \Amax := \exp (\Bbb R E). 
$$
For example,
 $\Mmax_+$ is the Lorentz group and $\Mmax_+ \Nmax$ is the Poincar\'e group
 if $(p,q) = (2,4)$.
It is convenient to identify $\Bbb R^n$ with $\overline{\Nmax}$ by putting
$$
   \nbar{a} := \exp (\sum_{j=1}^n a_j \overline{N_j}) \in \overline{\Nmax}
\qquad
 \text{for $a = (a_1, \dots, a_n) \in \Bbb R^n$.}
\tag \num.3
$$
The geometric point here will be the following:
\proclaim{Observation~\num}
{{\rm 1)}}
 On $S^{p-1} \times S^{q-1}$,
 $G$ acts conformally, while $K$ isometrically.
\newline{{\rm 2)}}
On $\Bbb R^{p-1,q-1}$,
 $G$ acts meromorphically and conformally,
 while the motion group $\Mmax_+ \overline{\Nmax}$  isometrically.
\endproclaim

\def \sc{6}
\sec{}
Next,
let us consider the pseudo-Riemannian manifold 
 $M = S^{p-1} \times S^{q-1}$.
It follows from (\ch.3.3) and (\ch.4.3) that we can define a representation
 $\varpi_{\lambda, M}$ of $G$ on $C^\infty(M)$ by
$$
   (\varpi_{\lambda, M} (h^{-1}) f) (v)
    := \nu(h \cdot v)^{-\lambda}  f(L_h v),  
$$

The Yamabe operator on $M$ is of the form:
$$
\tilLap{M} = \Delta_{S^{p-1}} - \Delta_{S^{q-1}} 
                - (\frac {p-2} 2)^2 + (\frac {q-2} 2)^2
 = \tilLap{S^{p-1}} - \tilLap{S^{q-1}}. 
$$
Applying Lemma~2.2, we obtain a  representation
 of the conformal group $G = O(p,q)$,
 denoted by $(\varpi^{p,q}, \Vpq)$, as a subrepresentation of
 $\varpi_{\frac{p+q-4}2, M}$:
$$
\alignat1
 &\Vpq := \Ker \tilLap{M}=
\set{f \in C^\infty(M)}{\tilLap{M} f =0},
\\
    &(\varpi^{p,q}(h^{-1})f) (v) := \nu(h \cdot v)^{-\frac {p+q-4} 2}
                                     f(L_h v),  
\quad \text{ for $h \in  G, v \in M, f \in \Vpq$}.  
\endalignat
$$

The restriction of $\varpi^{p,q}$
 from the conformal group to the isometry group
 gives a useful knowledge on the representation $\varpi^{p,q}$.
For this, we recall the classical theory of spherical harmonics,
 which is a generalization of Fourier series for $S^1$.
For $p \ge 2$ and $k \in \Bbb N$,
 we define the space of spherical harmonics of degree $k$ by
$$
\alignat1
     \spr{k}{p}
    &= \set{f \in C^{\infty}(S^{p-1})}
          {\Delta_{S^{p-1}} f = - k (k + p -2)f },
\tag \num.1
\\
    &= \set{f \in C^{\infty}(S^{p-1})}
          {\tilLap{S^{p-1}} f
     = \left(\frac{1}{4}- (k + \frac{p -2}2)^2\right) f }.
\endalignat
$$
Then $O(p)$ acts irreducibly on $\spr{k}{p}$ and
the algebraic direct sum
$
\bigoplus_{k=0}^\infty \spr{k}{p}
$
 is dense in $C^\infty(S^{p-1})$.
We note that $\spr{k}{2}\neq \{0\}$ only if $k = 0$ or $1$.

Now, we review a basic property of this representation $(\varpi^{p,q}, \Vpq)$
 on $M = S^{p-1} \times S^{q-1}$:
\proclaim{Lemma~\num\ {\rm (see \xbz; \xkorsI, \S 3)}} %
Assume $p, q \ge 2$, $p+q \in 2 \Bbb N$ and $(p,q) \neq (2,2)$.
\newline
{\rm 1)} 
 $(\varpi^{p,q},\Vpq)$ is an infinite dimensional irreducible
 representation of $G$.  
\newline
{\rm 2)}
 ($K$-type formula)\enspace
$\Vpq$ contains the algebraic direct sum
$$
 \bigoplus 
                     \Sb a, b \in \Bbb N \\ a + \frac{p}2 = b + \frac{q}2
                     \endSb
                      \spr{a}{p} \otimes \spr{b}{q} 
\tag \num.2
$$
 as a dense subspace with respect to the Fr\'echet topology on $C^\infty(M)$.
\newline
{\rm 3)}\enspace
 $G$ preserves the norm on $\Vpq$ defined by
$$
   \| F \|_M^2 := \|(\frac14-\tilLap{S^{p-1}})^{\frac14}F\|^2_{L^2(M)}
 =\sum_{a \ge \max(0,\frac{p-q}2)}
 (a + \frac{q-2}2) \|F_{a,b}\|^2_{L^2(M)}.
$$
if $F = \sum_{a} F_{a,b} \in \Vpq$ 
 with $F_{a,b} \in \spr{a}{p}\otimes \spr{b}{q}$ and $b = a + \frac{p-q}2$.
Here, 
$(\frac14-\tilLap{S^{p-1}})^{\frac14}$ is 
 a pseudo-differential operator on $M$,
 which is equal to
 $(\frac14-\tilLap{S^{q-1}})^{\frac14}$ on $\Ker \tilLap{M}$.
\endproclaim

We write $(\ , \ )_M$ for the corresponding inner product.
We denote by $\overline\Vpq$ the Hilbert completion of $\Vpq$,
 on which $G$ acts as an irreducible unitary representation of $G$.
We shall use the same notation $\varpi^{p,q}$ to denote this unitary
 representation.

If $p \ge q$ then $\Vpq$ contains the $K$-type
 of the form $\boldkey 1 \boxtimes \spr{\frac{p-q}2}{q}$.
This $K$-type is called a {\bf minimal $K$-type} in the sense of Vogan,
 namely,
 its highest weight 
 (with respect to a fixed positive root system of $\frak k_0$)
 attains the minimum distance
 from
 the sum of negative roots of $\frak k_0$
 among all highest weights of $K$-types occurring in $\varpi^{p,q}$.
Likewise for $p < q$.

\remark{Remark}
1)\enspace
If $p+q \ge 8$, $\varpi^{p,q}$ is called the {\bf minimal representation}
 in the representation theory of semisimple Lie groups,
 in the sense that the annihilator is the Joseph ideal.
\newline
2)\enspace
The formula (\num.2) is regarded as a branching law from
 the conformal group $G$ to the isometry subgroup $K$
 of the pseudo-Riemannian manifold $M = S^{p-1} \times S^{q-1}$
 (see Observation~2.5).
In \xkorsII, we generalized this branching law
 with respect to a non-compact reductive subgroup
 and proved the Parseval-Plancherel formula,
 in the framework of discretely decomposable restrictions \xkdecomp.
\endremark

\def\sc{7}
\sec{}
Let us consider the flat pseudo-Riemannian manifold $\Bbb R^{p-1,q-1}$.
The Yamabe operator on $\Bbb R^{p-1, q-1}$ is of the form:
$$
  \fLap \equiv \square_z := \frac{\partial^2}{{\partial z_1}^2} + \dots +
    \frac{\partial^2}{{\partial z_{p-1}}^2}
    - \frac{\partial^2}{{\partial z_p}^2} - \dots -
    \frac{\partial^2}{{\partial z_{p+q-2}}^2},
$$
 because the scalar curvature on $\Bbb R^{p-1,q-1}$ vanishes.
Since $G = O(p,q)$ acts on $\Bbb R^{p-1,q-1}$ 
 as a (meromorphic) conformal transform by (\ch.4.4),
 we obtain
 a \lq representation\rq\ with parameter $\lambda \in \Bbb C$ as in (\ch.2.3):
$$
 \varpi_{\lambda, \epsilon, \Bbb R^{n}}(g^{-1}) f(z)
  = |\mu(g \iota(z))|^{-\lambda} \chi_\epsilon(\sgn( \mu(g \iota(z))))
 f(L_g z),
 \quad (g \in G).
\tag \num.1
$$
Here,
 for $\epsilon = \pm 1$,
 we put
$$
\chi_\epsilon \: \Bbb R^\times \to \{\pm 1\}
$$ 
 by
 $\chi_{1} \equiv 1$ and $\chi_{-1} = \sgn$.
We may write 
$\varpi_{\lambda, \epsilon, \Bbb R^{p-1,q-1}}$ for 
 $\varpi_{\lambda, \epsilon, \Bbb R^{n}}$ 
 if we emphasize a view point of conformal geometry on
 the flat space $\Bbb R^{p-1,q-1}$.

We note that $C^\infty(M)$ is not stable by 
 $\varpi_{\lambda, \epsilon, \Bbb R^{p-1, q-1}}(g^{-1})$
 because $L_g$ is meromorphic.
To make (\num.1) a representation,
 we need to consider suitable class of 
 functions controlled at infinity.
One method for this is to use a conformal compactification
$$
  \Bbb R^{p-1, q-1} \hookrightarrow (S^{p-1} \times S^{q-1})/\sim \Bbb Z_2,
$$
 and to take a twisted pull-back $\Psi^*_\lambda$
 from $C^\infty(M)$ by a conformal map $\Psi$.
 This method is easy, and we shall explain it soon in \S \ch.8 and \S \ch.9.
The other is to find an inner product for specific parameter $\lambda$
 so that $G$ acts as a continuous unitary representation on the Hilbert space.
This is particularly non-trivial for a subrepresentation, 
 and we shall consider it for $\Ker \fLap$ in \S 6.

Before taking a suitable class of functions,
 we first write a more explicit form of (\num.1).
First, 
 we note that the maximal parabolic subgroup
$$
  \overline{\Pmax} := \Amax \Mmax \overline{\Nmax} 
 = (\Bbb R_+^\times \times O(p-1,q-1) \times \Bbb Z_2) \ltimes \Bbb R^n
$$
 acts transitively on the manifold $\iota(\Bbb R^n)$
 as affine transformations.
Furthermore, $\Mmax \overline{\Nmax}$ acts on $\iota(\Bbb R^n)$
 as isometries (see Observation~2.5).
Correspondingly,
 the representation
 $\varpi_{\lambda, \epsilon} 
 \equiv \varpi_{\lambda, \epsilon,\Bbb R^{n}}$ given in (\num.1)
 has a simple form
 when restricted
 to the subgroup $\overline\Pmax$:
$$
\alignat2
   (\varpi_{\lambda, \epsilon}(m) f)(z)
   &= f(m^{-1} z)
   && \qquad (m \in \Mmax_+),
\tag \num.2)(a
\\
 (\varpi_{\lambda, \epsilon}(m_0) f)(z)
   &= \epsilon f(z),
   && 
\\
   (\varpi_{\lambda, \epsilon} (e^{t E}) f)(z)
   &= e^{\lambda t} f (e^t z)
   && \qquad (t  \in \Bbb R),
\tag \num.2)(b
\\
   (\varpi_{\lambda, \epsilon} (\nbar{a}) f)(z)
   &= f(z-2 a)
   && \qquad (a \in \Bbb R^n).
\tag \num.2)(c
\endalignat
$$

Second,
 we write an explicit formula of the differential action of (\num.1).
We define a linear map
$$
   \omega \: \frak{g}_0 \to C^\infty(\Bbb R^n)
$$
by the Lie derivative of the conformal factor
$
   \Omega(h, z) := \mu(h \cdot \iota(z))^{-1}
$
 (see (2.4.4)).
For $Y =(Y_{i, j})_{0 \le i, j \le n+1} \in \frak{g}_0$
 and $z \in \Bbb R^n$, we have
$$
\omega(Y)_z 
:=  \frac{d}{d t}|_{t=0} \Omega(e^{t Y}, z) 
=
- Y_{0, n+1}
- \frac{1}2 \sum_{j=1}^{n} (Y_{0,j}+ Y_{n+1, j}) z_j.
\tag \num.3
$$

We write the Euler vector field on $\Bbb R^{n}$ as
$$
  E_z = \sum_{j=1}^{n} z_j \frac{\partial}{\partial z_j}.
\tag \num.4
$$
Then the differential
$
d \varpi_\lambda \:  \frak g_0 \to \operatorname{End}(C^\infty(\Bbb R^n))
$
 is given by
$$
\multline
d \varpi_\lambda(Y) 
=
-\lambda \omega(Y) -\omega(Y) E_z
\\
-
\sum_{i=1}^{n}
\left\{
  (Y_{i,0} + Y_{i, n+1})
 + \frac{|z'|^2 - |z''|^2}4 (-Y_{i,0} + Y_{i, n+1})
 +\sum_{j=1}^{n} Y_{i, j} z_j
\right\}
\frac{\partial}{\partial z_j}
\endmultline
\tag \num.5
$$
for $Y = (Y_{i, j})_{0 \le i, j  \le n+1} \in \frak g_0$
 and $z \in \Bbb R^n$.
In particular, we have
$$
d \varpi_\lambda ({N}_j) 
= 
- \lambda \varepsilon_j z_j 
- \varepsilon_j z_j E_z 
                     + \frac12 (|z'|^2 - |z''|^2)
 \frac{\partial}{\partial z_j}, 
\qquad (1 \le j \le n).
\tag \num.6
$$

\def\sc{8}
\sec{}{}
We recall $M = S^{p-1} \times S^{q-1}$.
This subsection relates the representation $\varpi_{\lambda, M}$ and 
 $\varpi_{\lambda, \Bbb R^{p-1,q-1}}$
 by the stereographic projection
 $\Psi^{-1} \: M \to \Bbb R^{p-1, q-1}$ defined below.

We set a positive valued function $\tau \: \Bbb R^n \to \Bbb R$ by
$$
\alignat1
\tau (z) \equiv \tau(z', z'') 
&
:= \nu \circ \iota(z)
\\
&=
 \left( (1 - \frac{|z'|^2 - |z''|^2}4)^2 + |z'|^2 \right)^{\frac12}
\\ 
&
=
 \left( (1 + \frac{|z'|^2 - |z''|^2}4)^2 + |z''|^2 \right)^{\frac12}
\\
&
=
 \left( 1 + (\frac{|z'| + |z''| }2)^2 \right)^{\frac12}
 \left( 1 + (\frac{|z'| - |z''| }2)^2 \right)^{\frac12}.
\tag \num.1
\endalignat
$$

We define an injective diffeomorphism
 as a composition of
  $\iota\: \Bbb R^{p-1,q-1} \hookrightarrow \Xi$ (see (2.3.4))
  and $\Phi \: \Xi \to M$ (see (2.3.5)):
$$
\Psi: \Bbb R^{p-1, q-1} \to M,
 \quad z \mapsto \tau(z)^{-1} \iota(z).
\tag \num.2
$$
The image of $\Psi$ is 
$$
   M_+ :=
 \set{u = (u_0, u', u'', u_{n+1}) 
\in M=S^{p-1} \times S^{q-1}}{u_0 + u_{n+1} > 0}.
\tag \num.3
$$
Then,
 $\Psi$ is a conformal map (see \xkorsI, Lemma 3.3, for example) such that
$$
    \Psi^* g_M = \tau(z)^{-2} g_{\Bbb R^{p-1,q-1}}.
\tag \num.4
$$
 The inverse of $\Psi \: \Bbb R^{p-1,q-1} \to M_+$ is given by
$$
\Psi^{-1} (u_0, u', u'', u_{n+1})
  = (\frac{u_0+u_{n+1}}2)^{-1} (u', u'')
  = {\mu(u)}^{-1} (u', u'')
\tag \num.5
$$
$\Psi^{-1}$ is nothing but a stereographic projection if $q = 1$.
We note that $\Psi$ induces a conformal compactification of
 the flat space $\Bbb R^{p-1, q-1}$:
$$
  \Bbb R^{p-1, q-1} \hookrightarrow (S^{p-1} \times S^{q-1})/\sim \Bbb Z_2.
$$
Here $\sim \Bbb Z_2$ denotes the equivalence relation in 
 the direct product space $S^{p-1} \times S^{q-1}$
 defined by $u \sim -u$.

As in (\ch.2.1),
 we define the twisted pull-back by
$$
\Psi^*_\lambda: C^\infty(M) \to C^\infty(\Bbb R^n),
\quad
 F \mapsto 
\tau(z)^{-\lambda} F(\Psi(z)).
\tag \num.6
$$
Let
$$
 C^\infty(M)_\epsilon :=
 \set{f \in C^\infty(M)}{f(- u) = \epsilon f(u), \ \text{ for any } u\in M}.
$$
Then $\Psi^*_\lambda|_{C^\infty(M)_\epsilon}$ is injective.
The inverse map is given by
$$
\multline
(\Psi_{\lambda, \epsilon}^*)^{-1} :
 \Psi^*_\lambda(C^\infty(M)_\epsilon) \to C^\infty(M)_\epsilon, 
\\
 f \mapsto
 \cases
  |\frac{u_0+u_{n+1}}2|^{-\lambda} f(\Psi^{-1}(u)) & (u \in M_+) 
 \\ 
 \epsilon
 |\frac{u_0+u_{n+1}}2|^{-\lambda} f(\Psi^{-1}(-u)) & (u \in M_-).
 \endcases
\endmultline
\tag \num.7
$$
We note that
$(\Psi_{\lambda, \epsilon}^*)^{-1} f$ makes sense for
 $f \in C_0^\infty(\Bbb R^n)$,
 since we have 
 $$C_0^\infty(\Bbb R^n) \subset \Psi^*_\lambda(C^\infty(M)_\epsilon).$$

Now, the representation 
 $\varpi_{\lambda,\epsilon, \Bbb R^{n}}$ is well-defined
 on the following representation space: $\Psi_\lambda^*(C^\infty(M))$,
 a subspace of $C^\infty(\Bbb R^{n})$
 through $\varpi_{\lambda, M}$.

Then, by (\ch.2.2) (see \xkorsI, Proposition~2.6), we have:
\proclaim{Lemma~\num}
$
\tilPsi(\Vpq) \subset \Ker \fLap,
$
where $\Vpq = \Ker \tilLap{M}$.
\endproclaim

\def\sc{9}
\sec{}
In the terminology of representation theory of semisimple Lie groups,
 $\Psi_\lambda^*$ is a $G$-intertwining operator from
 the $K$-{\it picture} $(\varpi_{\lambda, M}, C^\infty(M)_\epsilon)$
 to the $N$-{\it picture}
 $(\varpi_{\lambda, \epsilon,  \Bbb R^n}, \Psi_\lambda^*(C^\infty(M)))$.
To see this in an elementary way, we argue as follows:
For $\nu \in \Bbb C$,
 we denote by the space
$$
  S^\nu(\Xi) := \set{h \in C^\infty(\Xi)}{
                  h(t \xi) = t^\nu h(\xi), \ 
 \text{ for any } \xi \in \Xi, t > 0},
\tag \num.1
$$
 of smooth functions on $\Xi$ of homogeneous degree $\nu$.
Then $G$ acts on $S^\nu(\Xi)$ by left translations.
Furthermore,
 for $\epsilon = \pm 1$,
 we put
$$
 S^{\nu, \epsilon}(\Xi) := 
 \set{h \in S^\nu(\Xi)}{h(- \xi) = \epsilon h(\xi), \
 \text{ for any } \xi \in \Xi}.
\tag \num.2
$$
Then
 we have a direct sum decomposition
$$
  S^\nu(\Xi) = S^{\nu, 1}(\Xi) + S^{\nu,-1}(\Xi),
$$
 on which $G$ acts by left translations, respectively.
Then $S^{\nu, \epsilon}(\Xi)$ corresponds to the degenerate principal
 series representation (see \xkorsI\ for notation):
$$
 \prince{C^\infty}{\lambda}{\epsilon}
 \simeq S^{-\lambda -\frac{n}2, \epsilon}(\Xi),
\tag \num.3
$$
where $\Pmax = \Mmax \Amax \Nmax$.
\proclaim{Lemma~\num}
{\rm 1)}
The restriction $S^{-\lambda, \epsilon}(\Xi) \to C^\infty(M)_\epsilon,
 h \mapsto h|_M$
 induces 
 the isomorphism of $G$-modules between $S^{-\lambda, \epsilon}(\Xi)$ 
 and $(\varpi_\lambda, C^\infty(M)_\epsilon)$ for any $\lambda \in \Bbb C$.
\newline
{\rm 2)}\enspace
The restriction $S^{-\lambda, \epsilon}(\Xi) \to C^\infty(\Bbb R^n),
 h \mapsto h|_{\Bbb R^n}$
 induces 
 the isomorphism of $G$-modules between $S^{-\lambda, \epsilon}(\Xi)$ 
 and $(\varpi_{\lambda, \epsilon, \Bbb R^n}, \Psi_\lambda^*(C^\infty(M)_\epsilon))$
 for any $\lambda \in \Bbb C$.
\endproclaim
\demo{Proof}
See \xkorsI, Lemma~3.7.1 for (1).
 (2) follows from the commutative diagram:
$$
\alignat3
&&&\hphantom{M}  S^{-\lambda, \epsilon}(\Xi) &&
\\
&&r_1\swarrow&& \searrow r_2&
\tag \num.4
\\
&&C^\infty(M)_\epsilon
&\ \quad \ \overset{\Psi_\lambda^*}\to\longrightarrow && C^\infty(\Bbb R^n),
\endalignat
$$
because $r_1$ is bijective
and
 $r_2$ is injective.
\qed
\enddemo
\def\sc{10}
\sec{}
A natural bilinear form
$
\langle \ , \ \rangle :
S^{-\lambda-\frac{n}2}(\Xi) \times S^{\lambda-\frac{n}2}(\Xi) \to \Bbb C
$
 is
defined by
$$
\alignat1
  \langle h_1, h_2 \rangle 
  &:=\int_M h_1(b) h_2(b) \ d b
\tag \num.1
\\
  &= 2 \int_{\Bbb R^n} h_1(\iota(z)) h_2(\iota(z)) \ d z 
      \qquad \text{(see (\ch.3.4))}.
\tag \num.2
\endalignat
$$
Here, $d b$ is the Riemannian measure on $M = S^{p-1} \times S^{q-1}$.
The second equation follows from 
$
   (h_1 h_2)(\iota(z)) =
   \tau(z)^{-n}
   (h_1 h_2)(\Psi(z))
$
 and the Jacobian for $\Psi \: \Bbb R^n \to M$ is given by $\tau(z)^{-n}$
 (see (\ch.8.4)).
Then $\langle \ , \ \rangle$ is $K$-invariant and $\overline{\Nmax}$-invariant
 from (\num.1) and (\num.2),
 and thus $G$-invariant since $G$ is generated by $K$ and $\overline{\Nmax}$.
\def\ch{3}
\head
\S \ch.  
   Square integrable functions on the cone
\endhead
\def\sc{1}
\sec{}
In this section we shall study the irreducible unitary representation of
the motion group
 $\Mmax_+ \overline\Nmax \simeq O(p-1, q-1) \ltimes \Bbb R^{p+q-2}$ and the
maximal parabolic subgroup ${\overline\Pmax} = \Mmax \Amax \overline\Nmax$ 
on the space of solutions to
our ultrahyperbolic equation $\fLap f = 0$. 
This is a standard induced representation by the Mackey machine, and
will be later extended to
the minimal representation of $G = O(p,q)$ (see Theorem~4.9~(3)).  

\def\sc{2}
\sec{}
In the flat picture $\Bbb R^{p-1. q-1}$,
 our minimal representation $\Vpq$ of $O(p,q)$ can be realized in a subspace of 
 $\Ker \fLap$ (see Lemma~2.8).
We shall study the representation space by means of the Fourier transform.

We normalize the Fourier transform on $\Cal S(\Bbb R^n)$ by
$$
   (\Cal F f)(\zeta) %
       = \int_{\Bbb R^n} f(z) e^{\sqrt{-1} (z_1 \zeta_1 + \dots + z_n \zeta_n)}
                            d z_1 \cdots d z_n,        
$$
 and extends it to $\Cal S'(\Bbb R^n)$,
 the space of the Schwartz distributions.

By composing the following two injective maps
$$
 C^\infty(M)_\epsilon  \overset{\tilPsi}\to\longrightarrow 
              C^\infty(\Bbb R^n) \cap \Cal S'(\Bbb R^n)
              \overset{\Cal F}\to\longrightarrow 
              \ \Cal S'(\Bbb R^n),
$$
 we define a representation of $G$ and $\frak g$
 on the image
 $\Cal F \tilPsi (C^\infty(M))$,
 denoted by
 $\widehat\varpi_{\lambda, \epsilon}
 \equiv \widehat\varpi_{\lambda, \epsilon, \Bbb R^{n}}$,
 so that
 $\Cal F \circ \tilPsi$
 is a bijective $G$-intertwining operator
 from the representation space $(\varpi_{\lambda, M}, C^\infty(M)_\epsilon)$
 to $(\widehat\varpi_{\lambda, \epsilon}, \Cal F \tilPsi(C^\infty(M)_\epsilon))$.
Then, 
 it follows from (2.7.2) that
 the representation
 $\widehat\varpi_{\lambda, \epsilon}$
 has a simple form
 when restricted to the subgroup $\overline\Pmax=\Mmax \Amax \overline\Nmax$:
$$
\alignat2
   (\widehat\varpi_{\lambda, \epsilon}(m) h)(\zeta)
   &= h(\trans{m} \zeta)
   && \qquad (m \in \Mmax_+),
\tag \num.1)(a
\\
   (\widehat\varpi_{\lambda, \epsilon}(m_0) h)(\zeta)
   &= \epsilon h(\zeta),
   && 
\\
   (\widehat\varpi_{\lambda, \epsilon}(e^{t E}) h)(\zeta)
   &= e^{(\lambda -n) t} h(e^{-t} \zeta)
   && \qquad (t \in \Bbb R),
\tag \num.1)(b
\\
   (\widehat\varpi_{\lambda, \epsilon}(\nbar{a}) h)(\zeta)
   &= e^{2\sqrt{-1} (a_1 \zeta_1 + \dots + a_n \zeta_n)} h(\zeta)
   && \qquad (a \in \Bbb R^n).
\tag \num.1)(c
\endalignat
$$
We remark that in the above formula,
 we regarded  $h$ as a function.
The action of $\Amax$ on the space of distributions is slightly different
 by the contribution of the measure $d \zeta$:
$$
   (\widehat\varpi_{\lambda, \epsilon}(e^{t E}) \phi)(\zeta)
  = e^{\lambda t} \phi(e^{-t} \zeta)
   \qquad (t \in \Bbb R),
\tag \num.1)(b$'$
$$
if we write $\phi(\zeta) = h(\zeta) d \zeta \in \Cal S'(\Bbb R^n)$.

The differential representation $d \widehat\varpi_{\lambda, \epsilon}$
 of $\frak g_0$
 on $\Cal F \tilPsi (C^\infty(M))$ is given by the following lemma:
\proclaim{Lemma~\num}
We recall that $E_\zeta$ is the Euler operator (see (2.7.4)).
With notation in (2.5.2), we have
$$
\alignat2
d \widehat\varpi_{\lambda, \epsilon} (\overline{N}_j) 
&= 2 \sqrt{-1} \zeta_j
\qquad && (1 \le j \le n),
\\
d \widehat\varpi_{\lambda, \epsilon} ({N}_j) 
&= \sqrt{-1} 
\left( \lambda \varepsilon_j \frac{\partial}{\partial \zeta_j} 
- E_\zeta \varepsilon_j \frac{\partial}{\partial \zeta_j} 
                     + \frac12 \zeta_j \square_\zeta\right)
\qquad && (1 \le j \le n),
\\
d \widehat\varpi_{\lambda, \epsilon} (E) &= \lambda -n-E_\zeta.
\endalignat
$$
\endproclaim
\demo{Proof}
Lemma follows from the correspondence under the Fourier transform
 $\frac{\partial}{\partial z_j} \leftrightarrow -\sqrt{-1} \zeta_j$,
 $z_j \leftrightarrow -\sqrt{-1}\frac{\partial}{\partial \zeta_j}$,
 and therefore from
 $E_z \leftrightarrow -n-E_\zeta$,
 $P(z) \leftrightarrow -\square_\zeta$,
 where $P(z)
 := {z_1}^2 + \dots + {z_{p-1}}^2 - {z_p}^2 -\dots-{z_{p+q-2}}^2$.
\qed
\enddemo
\remark{Remark}
{\rm  1)}\enspace
We note that 
 $d \widehat\varpi_{\lambda, \epsilon}$ is independent of
 the signature $\epsilon = \pm 1$.
\newline
{\rm  2)}\enspace
In Theorem~4.9,
 we shall find that $L^2(C)$,
 the Hilbert space of square integrable functions on the cone $C$
 in $\Bbb R^n$,
 is a $G$-invariant subspace of the Schwartz distributions $\Cal S'(\Bbb R^n)$.
Then, the action of the Lie algebra $\frak g$
 can be written in terms of differential operators along the cone $C$
 at most of second order.
\endremark

\def\sc{3}
\sec{}
We define a quadratic form $Q$ on $\Bbb R^n$ ($\simeq (\Bbb R^n)^*$) 
 as the dual of $P(z)$ on $\Bbb R^n$ by 
$$
   Q(\zeta) := {\zeta_1}^2 + \dots + {\zeta_{p-1}}^2
                    - {\zeta_{p}}^2 - \dots - {\zeta_{p+q-2}}^2
\qquad
\text{ for $\zeta \in \Bbb R^n = \Bbb R^{p+q-2}$}
\tag \num.1
$$
 and define a closed cone by
$$
    C := \set{\zeta \in \Bbb R^n}{Q(\zeta) =0}.
\tag \num.2
$$
It follows from Lemma~2.8 that the support of 
 the distribution $\Cal F \tilPsi F$
 is contained in the cone $C$, for any $F \in \Vpq$.
Surprisingly,
 $\Cal F \tilPsi F$
 becomes square integrable on $C$ (see Theorem 4.9).
As a preparation for the proof, we study a natural action on $L^2(C)$
 of a parabolic subgroup $\overline{\Pmax}$ in this subsection.

We take a differential $(n-1)$-form $d \mu$ on $C$ such that
$$
 d Q \wedge d \mu = d \zeta_1 \wedge \dots \wedge d \zeta_n.
$$
Then the restriction $d \mu$ to the cone $C$ defines a canonical measure
 (we use the same notation $d \mu$).
Using polar coordinates on $C$:
 $\zeta = (s \omega, s \omega')$ with $s > 0$, $\omega \in S^{p-2}$,
 $\omega' \in S^{q-2}$,
 we write down the canonical measure $d \mu$ on $C$ explicitly by
$$
   \int_C \phi d \mu 
  =
\frac{1}{2}
 \int_0^\infty \int_{S^{p-2}} \int_{S^{q-2}} \phi(s \omega, s \omega')
          s^{n-3}\ d s \ d \omega \ d \omega'
\tag \num.3
$$
 for a test function $\phi$.
If $n > 2$, that is if  $p+q>4$, then the measure $d \mu$ defines 
 a Schwartz  distribution on $\Bbb R^n$,
 which is equal to the Delta function $\delta (Q)$
 supported on the cone $C$ (see \xgs).

For a measurable function $\phi$ on $C$,
 we define a norm of $\phi$ by
$$
   \| \phi\|^2 :=
   \int_C |\phi|^2 d \mu 
\tag \num.4
$$
 and denote by $L^2(C) \equiv L^2(C, d \mu)$
 the Hilbert space of square integrable functions.

By the Mackey theory,
 we can define a natural representation 
 $\pi$ of
 the maximal parabolic subgroup $\overline\Pmax=\Amax \Mmax \overline{\Nmax}$ on $L^2(C)$
 if $p+q \in 2 \Bbb Z$,
  by
$$
\alignat2
   (\pi(e^{t E}) \psi)(\zeta)
   &:= e^{-\frac{n-2}2 t} \psi(e^{-t} \zeta),
   && \qquad (t \in \Bbb R)
\tag \num.5)(a
\\
   (\pi(m) \psi)(\zeta)
   &:= \psi(\trans{m} \zeta)
   && \qquad (m \in \Mmax_+),
\tag \num.5)(b
\\
 (\pi(m_0) \psi)(\zeta)
   &:= (-1)^{\frac{p-q}2}\psi(\zeta)
   && 
\\
   (\pi(\nbar{a}) \psi)(\zeta)
   &:= e^{2\sqrt{-1} (a_1 \zeta_1 + \dots + a_n \zeta_n)} \psi(\zeta)
   && \qquad (a \in \Bbb R^n).
\tag \num.5)(c
\endalignat
$$
\proclaim{Proposition~\num}
{\rm 1)}\enspace
The representation $(\pi, L^2(C))$ of $\overline{\Pmax}$ is unitary.
\newline
{\rm 2)}\enspace
The representation $\pi$ is still irreducible when restricted to 
 the motion group
 $\Mmax_+ \overline{\Nmax} \simeq O(p-1,q-1)\ltimes \Bbb R^{p+q-2}$.
In particular, it is irreducible as a $\overline{\Pmax}$-module.
\endproclaim
\demo{Proof}
(1) is straightforward from the definitions (\num.4) and (\num.5).
\newline
Let us prove (2).
It follows from (\num.5)(c) that
 any $\overline{\Nmax}$-invariant closed subspace of $L^2(C)$ is of the form
 $L^2(C')$ where $C'$ is a measurable subset of $C$.
As $\Mmax_+$ acts transitively on $C$,
 $L^2(C')$ is $\Mmax_+$-invariant only if
 the measure of $C'$ is either null or conull.
Thus, $L^2(C')$ equals either $\{0\}$ or $L^2(C)$.
Therefore,
 the unitary representation $L^2(C)$ is irreducible 
 as an $\Mmax_+ \overline{\Nmax}$-module.
\qed
\enddemo

\def\sc{4}
\sec{}
It is not clear a priori if $(\pi, L^2(C))$
 extends from $\overline{\Pmax}$ to $G$.
We shall prove in Theorem~4.9 that 
 if $p, q \ge 2$ and $n (=p+q-2)>2$ then
 $\pi$ extends to $G$
 as an irreducible unitary representation 
 through an injective map
$
    T \: L^2(C) \hookrightarrow \Cal S'(\Bbb R^n), 
$
defined as follows:

By using the Cauchy-Schwarz inequality,
 we see the following map
$$
 T(\psi) \: \Cal S(\Bbb R^n) \to \Bbb C,
 \quad     \varphi \mapsto
 \int_C \varphi \psi \ d\mu
$$
 is well-defined and continuous if $n >2$, 
 for each $\psi \in L^2(C)$.
Thus we have a natural map
$$
    T \: L^2(C) \longrightarrow \Cal S'(\Bbb R^n), 
 \quad \psi \mapsto \psi d \mu.
\tag \num.1
$$
Clearly, $T$ is injective.
We shall regard $T(L^2(C))$ as a Hilbert space such that $T$ is a unitary
 operator.
\proclaim{Lemma~\num}
{\rm 1)}\enspace
$T$ is a $\overline\Pmax$-intertwining operator
 from $(\pi, L^2(C))$
 to 
\newline
$(\widehat\varpi_{\frac{n-2}2, \epsilon}|_{\overline{\Pmax}}, \Cal S'(\Bbb R^n))$.
\newline
{\rm 2)} \enspace
 $(\widehat\varpi_{\frac{n-2}2, \epsilon}|_{\overline{\Pmax}}, T(L^2(C)))$
 is an irreducible unitary representation of $\overline{\Pmax}$.
It is still irreducible as an $\Mmax \overline{\Nmax}$-module.
\endproclaim
\demo{Proof}
(1) follows directly from the definitions (\ch.2.1) and (\ch.3.5).
(2) follows from (1) and Proposition~\ch.3~(2).
\qed
\enddemo
\def\ch{4}
\head \S \ch. Green function and inner product \endhead
\def\sc{1}
\sec{}
In this section,
 we shall construct solutions of the ultrahyperbolic equation
 $\fLap f = 0$
 by the integral transform given by convolution with the Green kernel.
Then, we shall show that the Green kernel coincides with
 a special value of the Knapp-Stein intertwining operator for a degenerate
 principal series.
This observations gives another expression of the inner product
 of the minimal representation of $O(p,q)$ by using the Green kernel
 (see Theorem~4.7),
 and also leads to a realization of the minimal representation on
 $L^2(C)$, the Hilbert space of square integrable functions on a cone $C$
 as will be discussed in \S 6.

We put
$$
 P(x) = 
{x_1}^2 + \dots + {x_{p-1}}^2
 - {x_p}^2 - \dots -{x_{n}}^2
$$
for $x \in \Bbb R^n = \Bbb R^{p+q-2}$.
\def\sc{2}
\sec{}
A distribution $E$ satisfying $\fLap E = \delta$ (Dirac's delta function)
 is called a {\it fundamental solution} of the ultra-hyperbolic Laplace
 operator $\fLap$.
Recall from \xgs, page 354,
 if $n$ is even and $n > 2$ then
$$
  E = 
\frac{- \Gamma(\frac{n}2-1) e^{\frac{\sqrt{-1} \pi (q-1)}2}}{4 \pi^{\frac{n}2}}
 (P(x) + \sqrt{-1} 0)^{1-\frac{n}2}
$$
 is a fundamental solution of $\fLap$,
 where $(P(x) + \sqrt{-1} 0)^\lambda$ stands for the limit of
the distribution $(P(x) + \sqrt{-1} R(x))^\lambda$ as a positive definite
 quadratic form $R(x)$ tends to $0$.
In view of the integral formula in \xgs, Chapter III, \S 2.6
$$
  \Cal F (P + \sqrt{-1} 0)^{1 - \frac{n}2}
 = \frac{4 \pi^{\frac{n}2} e^{-\frac{\sqrt{-1} \pi(q-1)}{2}}}{\Gamma(\frac{n}2-1)} (Q -\sqrt{-1}0)^{-1},
$$
we have readily the following formula for the Green function 
$E_0$ of $\fLap$:
\proclaim{Proposition~\num}
We define a distribution $E_0$ on $\Bbb R^n$ by 
$$
      E_0 =\frac{1}{2 \pi \sqrt{-1}} (-E + \overline{E}).
\tag \num.1
$$
Then
 $\fLap E_0  = 0$
 and its Fourier transform is given by
 
$$
   \Cal F{E_0} =
 \frac{1}{2\pi\sqrt{-1}}((Q-\sqrt{-1}0)^{-1} - (Q + \sqrt{-1}0)^{-1})
= \delta(Q).
\tag \num.2
$$
\endproclaim

In the Minkowski case, i.e. $q = 2$,
 such a formula has been known (\xabste),
since the so-called two-point functions in the quantum field theory 
for a zero mass field exactly corresponds to $(Q + \sqrt{-1}0)^{-1}$
for negative frequency. 
In this case $C$ naturally splits in
two components, a forward and a backward light cone, and functions
supported on the forward cone have Fourier transforms that extend to
holomorphic functions on the corresponding tube domain, thus
yielding a unitary highest weight representation of the connected group.
The reproducing kernel of this representation is the Fourier
transform of the measure on the forward cone, in analogy with
what happens in Proposition 4.2.   
   
\def\sc{3}
\sec{}
In order to give the integral expression of solutions $\fLap f = 0$,
 we define a convolution map by the Green kernel:
$$
   S \: C^\infty_0(\Bbb R^n) \to C^\infty(\Bbb R^n),
  \quad
  \varphi \mapsto E_0 * \varphi. 
\tag \num.1
$$
Fix $\nu \in \Bbb C$.
We consider the representation $\varpi_{\nu, \epsilon} \equiv
 \varpi_{\nu, \epsilon, \Bbb R^{n}}$
 of $G$ on a subspace of $C^\infty(\Bbb R^n)$ (see \S 2).
The restriction to $\overline\Pmax$ stabilizes $C^\infty_0(\Bbb R^n)$,
 as follows from (2.7.2).
\proclaim{Lemma~\num}
{\rm 1)}\enspace
$\Image S \subset \Ker \fLap$.
\newline
{\rm 2)}\enspace
$S$ is an intertwining operator of $\overline{\Pmax}$-modules
 between on one side $\varpi_{\nu+2, \epsilon}|_{\overline\Pmax}$ 
 and on the other side $\varpi_{\nu, \epsilon}|_{\overline\Pmax}$
 for any $\nu \in \Bbb C$ and $\epsilon = \pm 1$.
\endproclaim
\demo{Proof}
(1)\enspace
As $\fLap E_0 =0$,
 we have
 $\Image S \subset \Ker \fLap$.
\newline
The proof of (2) is direct from (2.7.2).
We illustrate it by the action of $e^{t E} \in \Amax$:
$$
\alignat1
S(\varpi_{\nu+2, \epsilon}(e^{t E}) \varphi)(y) 
& =
\int_{\Bbb R^n} E_0(y-z) e^{(\nu + 2)t} \varphi(e^t z) \ d z
\\
& =
\int_{\Bbb R^n} E_0(e^t y- e^t z) e^{(\nu + 2)t+(n-2)t} \varphi(e^t z) 
 \ d z
\\
& =
e^{\nu t} (S \varphi)(e^t y)
\\
& = \varpi_{\nu, \epsilon}(e^{t E}) (S \varphi)(y).
\endalignat
$$
This shows that $S$ intertwines the action of $\Amax$.
The case for the action of $\Mmax \overline{\Nmax}$ is similar and easier.
\qed
\enddemo
We shall see that $S$ extends to a $G$-intertwining operator for $\nu = \frac{n}2-1$
in Proposition~\ch.6.

\def\sc{4}
\sec{}
Recall the notation in \S 2.9.
Let
$$
[u,v] := u_0 v_0 + \dots + u_{p-1} v_{p-1}
 - u_p v_p - \dots -u_{n+1} v_{n+1}
$$
 for $u, v \in \Bbb R^{n+2} = \Bbb R^{p+q}$.
The Knapp-Stein intertwining operator
$$
A_{\lambda, \epsilon} \:
  S^{-\lambda - \frac{n}2, \epsilon}(\Xi)
 \to  S^{\lambda - \frac{n}2, \epsilon}(\Xi)
$$
 is given by the integral operator with kernel function 
$$
 \psi_{\lambda -\frac{n}2, \epsilon}(u_0-u_{n+1})
 = \psi_{\lambda -\frac{n}2, \epsilon}([u, \xi_0])
\tag \num.1
$$
on $\Xi$.
Here,
 we put $\xi_0 := \trans (1, 0, \dots, 0,1)\in \Xi$
 and
 $\psi_{\nu, \epsilon}$ is a distribution (or a hyperfunction) of one variable
 is defined by
$$
  \psi_{\nu, \epsilon}(y) := \frac{1}{\Gamma(\frac{2 \nu + 3 - \epsilon}{4})}
                              |y|^\nu \chi_\epsilon(\sgn y).
$$

Via the bijection $r_1$ in the commutative diagram (2.9.4),
$$
A_{\lambda, \epsilon}: C^\infty(M)_\epsilon \to C^\infty(M)_\epsilon
$$ 
 is written as
$$
 (A^{M}_{\lambda, \epsilon} f)(u)
\equiv (A_{\lambda, \epsilon} f)(u)
 := \int_M
 \psi_{\lambda-\frac{n}2, \epsilon}([u, v]) \ f(v) \ d v
  \quad
  (u \in M)
$$
in the compact picture $M \simeq S^{p-1} \times S^{q-1}$. 

The Gamma factor in the definition of
 $\psi_{\nu, \epsilon}(y)$ exactly cancels the poles of 
  the distribution $|y|^\nu \chi_\epsilon(\sgn y)$ of one variable $y$
 with meromorphic parameter $\nu$.
This means that
 the distribution $\psi_{\nu, \epsilon}([u,v])$ of multi-variables makes sense
 for any $\nu \in \Bbb C$
 when restricted to the open set
 $\set{(u,v) \in M \times M}{u \neq \pm v}$,
  where $[u, dv] + [v, d u] \neq 0$.
Then
 $\psi_{\nu, \epsilon}([u,v])$ 
 continues meromorphically as a distribution on $M\times M$
 with possible poles only at $\lambda = 0, -1, \dots$,
 whose residues are distributions supported on
 $\set{(u,v) \in M \times M}{u = \pm v}$.

In view of the normalization of our parameter (see Lemma~2.9),
 $A_{\lambda, \epsilon}$ is a $G$-intertwining operator
 from $\varpi_{\lambda+\frac{n}2}$
 to $\varpi_{-\lambda+\frac{n}2}$
 for $\lambda \neq 0, -1, -2, \dots$.
What we need is the case $\lambda =1$ and
we recall from \xkorsI, \S 3.9
 (basically since we know the composition series
 of the induced representations
 and the eigenvalue of $A_{\lambda, \epsilon}$ on each $K$-type):
\proclaim{Lemma~\num}
 Let $p \equiv q \mod 2$, $p, q \ge 2$ and $(p,q) \neq (2,2)$.
We put
$$
\alignat2
\epsilon &:= (-1)^{\frac{p-q}2}
&&= \cases 1 & p-q\equiv 0, 4 \mod 8,
                                       \\
                                -1 & p-q\equiv 2, 6 \mod 8. \endcases
\tag \num.2)(a
\\
\delta &:= (-1)^{[\frac{q-p}4]} 
&&= \cases 1 & p-q\equiv 0, 2 \mod 8,
                                       \\
                                -1 & p-q\equiv 4, 6 \mod 8. \endcases
\tag \num.2)(b
\endalignat
$$
{\rm 1)}\enspace
The image  $(\varpi_{\frac{n}2-1}, A_{1, \epsilon}(C^\infty(M)_\epsilon))$
 is a dense subrepresentation of $(\varpi^{p,q}, \Vpq)$,
  where we recall $\Vpq = \Ker \tilLap{M}$.

Let us define a Hermitian form $( \ , \ )_A$ on the same image
 $A_{1, \epsilon}(C^\infty(M)_\epsilon)$ by
$$
 (F_1, F_2)_A
:=
\delta  \langle A_{1, \epsilon} \varphi_1, \overline{\varphi_2} \rangle
=\delta \langle  \varphi_1, \overline{A_{1, \epsilon} \varphi_2}\rangle,
\tag \num.3
$$
for $F_i = A_{1, \epsilon} \varphi_i$, $\varphi_i \in C^\infty(M)_\epsilon$
 ($i=1, 2$),
 where $\langle \ , \ \rangle$ is the bilinear form given 
 as the integral over $M$ (see (2.10.1)).
\newline{\rm 2)}\enspace
The Hermitian form 
 $(\ , \ )_A$ is well-defined (namely, independent of the choice of $\varphi_i$)
 and $G$-invariant under the action of $\varpi_{\frac{n}2-1}$.
\newline{\rm{3)}\enspace
In comparison with the inner product $(\ , \ )_M$ given in Lemma~2.6, we have
$$
(F_1, F_2)_A = c_1 (F_1, f_2)_M
\qquad
\text{ for any } F_1, F_2 \in A_{1, \epsilon}(C^\infty(M)_\epsilon),
\tag \num.4
$$
where we put
$$
 c_1 := \frac{\Gamma(\frac{n-1-\epsilon}4)}{2^{\frac{n}2} \pi^{\frac{n+1}2}}.
$$
In particular, the Hermitian form $(\ , \ )_A$
 is positive definite
 and the completion of a pre-Hilbert space 
  $(A_{1, \epsilon}(C^\infty(M)_\epsilon), (\ , \ )_A)$
 coincides with the Hilbert space $\overline\Vpq$ given in \S 2.6.
\endproclaim
\def\sc{5}
\sec{}
In the flat picture $\Bbb R^{p-1,q-1}$,
 we have 
$$
     -\xi_0+\xi_{n+1} = |z'|^2-|z''|^2 = P(z)
$$
 for $\xi = \iota(z)$ by (2.3.4).
Then,
 via the injection $r_2$ in the commutative diagram (2.9.4),
 the Knapp-Stein intertwining operator is given by the convolution:
$$
 A^{\Bbb R^n}_{\lambda, \epsilon} \varphi
\equiv
 A_{\lambda, \epsilon} \varphi
 := 2 \psi_{\lambda-\frac{n}2, \epsilon}(P(z)) * \varphi
\tag \num.1
$$
 when restricted to
 $C^\infty_0(\Bbb R^n) \subset \Psi^*_{\frac{n+2}2}(C^\infty(M)_\epsilon)$.
Then we have
$$
 A^{\Bbb R^n}_{1, \epsilon} \circ \Psi^*_{\frac{n+2}2}
=
 \tilPsi \circ A^{M}_{1, \epsilon} .
\tag \num.2
$$
\proclaim{Lemma~\num}
Retain the setting of Lemma~\ch.4.
We put
$$
      h(y) := e^{\frac{\sqrt{-1}\pi(q-1)}{2}} (y+\sqrt{-1}0)^{1-\frac{n}2}.
\tag \num.3
$$
Then we have
$$
\psi_{1-\frac{n}2,\epsilon}(y)
 = \delta \frac{\Gamma(\frac{-1+n+\epsilon}{4})}{2\pi\sqrt{-1}}
 (h(y) - \overline{h(y)}).
\tag \num.4
$$
\endproclaim
\demo{Proof}
Solving
 $(y\pm\sqrt{-1} 0)^\nu = y_+^\nu + e^{\pm\sqrt{-1} \pi \nu} y_-^\nu$,
 we have
$$
  \psi_{\nu, \epsilon}(y) = 
\frac{(e^{-\sqrt{-1}\nu\pi} -\epsilon)(y+\sqrt{-1}0)^\nu
- (e^{\sqrt{-1}\nu\pi} -\epsilon)(y-\sqrt{-1}0)^\nu}{
-2 \sqrt{-1}\sin\pi\nu \Gamma(\frac{2 \nu + 3 - \epsilon}{4})}.
$$
Then, Lemma~\num\ follows from a residue calculation:
$$
    \lim_{\nu \to 1 - \frac{n}2}
\frac{e^{-\sqrt{-1}\nu\pi} -\epsilon}{
\sin\pi\nu \ \Gamma(\frac{2 \nu + 3 - \epsilon}{4})}
=
   \lim_{a\to 0}
\frac{e^{-\sqrt{-1}a\pi} -(-1)^q}{
\sin\pi a \ \Gamma(\frac{5-n-\epsilon}{4}+\frac{a}2)}
=\frac{-\delta e^{\frac{\sqrt{-1}\pi(q-1)}{2}} \Gamma(\frac{-1+n+\epsilon}{4})}{\pi}.
$$
\enddemo
\def\sc{6}
\sec{}
Our key observation is that the special value of the Knapp-Stein operator
 is given by the Green function up to a scalar constant.
In particular,
 the image of $C^\infty_0(\Bbb R^n)$ under the integral transform $S$
 satisfies not only the ultra-hyperbolic equation $\fLap f = 0$
 but also a certain decay condition at infinity
 that matches a conformal compactification
 of $\Bbb R^{p-1, q-1}$ (see \S 2.8).
\proclaim{Proposition~\num}
Retain the notation as in Proposition~\ch.2 and Lemma~\ch.4. 
We recall that $\delta, \epsilon = \pm 1$ are determined by $p-q \mod 8$
as in (4.4.2) and define the constant $c_2$ by
$$
c_2 := 
 \frac{4\delta \pi^{\frac{n}2}\Gamma(\frac{-1+n+\epsilon}4)}{\Gamma(\frac{n}{2}-1)}.
$$
Then we have:
$$
\alignat1
  \psi_{1-\frac{n}2, \epsilon}(P(z))
   &= c_2 E_0,
\tag \num.1
\\
   A^{\Bbb R^n}_{1, \epsilon}|_{C^\infty_0(\Bbb R^n)}
   &= 2 c_2 S.
\tag \num.2
\endalignat
$$
In particular,
 $S (C_0^\infty(\Bbb R^n)) \subset \tilPsi(\Vpq)$.
\endproclaim
\demo{Proof}
The first formula follows from Lemma~\ch.5 and from
 the definition of $E_0$ (see (\ch.2.1)).
The second formula then follows from the definition (\ch.5.1).
Then, in view of (\ch.5.2) and (\num.2),
 $S$ extends to a $G$-intertwining operator 
 between $(\varpi_{\frac{n+2}2}, \Psi_{\frac{n+2}2}^*(C^\infty(M)_\epsilon))$
 and $(\varpi_{\frac{n-2}2}, \tilPsi(C^\infty(M)_\epsilon))$.
Then we have
$$
 S (C_0^\infty(\Bbb R^n)) \subset
 A^{\Bbb R^n}_{1, \epsilon} (\Psi_{\frac{n+2}2}^*(C^\infty(M)_\epsilon))
 =\tilPsi ( A^{M}_{1, \epsilon}(C^\infty(M)_\epsilon))
 \subset 
 \tilPsi(\Vpq).
$$
 because
 $C^\infty_0(\Bbb R^n) \subset \Psi_{\frac{n+2}2}^*(C^\infty(M)_\epsilon)$.
\qed
\enddemo
\def\sc{7}
\sec{}
We are ready to introduce
 an inner product on $\Ker \fLap$ by using the Green function $E_0$
 (see (\ch.2.1)),
 so that our minimal representation is realized here as
 a unitary representation.

We define a Hermitian form on $S(C^\infty_0(\Bbb R^n))
 \subset \Ker \fLap$ by
$$
  (f_1, f_2)_N :=
 \int_{\Bbb R^n}\int_{\Bbb R^n} E_0(y-x) \varphi_1(x) \overline{\varphi_2(y)}
 \ d x \ d y
\quad
\tag \num.1
$$
for $f_i = S \varphi_i = E_0 * \varphi_i$ ($1 \le i \le 2$).
The right-hand side of (\num.1) does not depend on the choice of $\varphi_i$
 because of the formula
$$
   (f_1, f_2)_N = (f_1, \varphi_2)_{L^2(\Bbb R^n)}
                = (\varphi_1, f_2)_{L^2(\Bbb R^n)}.
$$
We recall from Proposition~\ch.6 and Lemma~2.8 the inclusive relations:
$$
S (C_0^\infty(\Bbb R^n)) \subset 
A^{\Bbb R^n}_{1, \epsilon}
 (\Psi_{\frac{n+2}2}^*(C^\infty(M)_\epsilon)) \subset 
\tilPsi(\Vpq)
\subset \Ker \fLap \subset C^\infty(\Bbb R^n).
\tag \num.2
$$
\proclaim{Theorem~\num}
 Let $p \equiv q \mod 2$, $p, q \ge 2$ and $(p,q) \neq (2,2)$.
Recall $n = p+q-2$ and $\epsilon=(-1)^{\frac{p-q}2}$.
\newline
{\rm 1)}\enspace
The Hermitian form $(\ , \ )_N$ is positive definite on 
 $S(C^\infty_0(\Bbb R^n))$.
\newline
{\rm 2)}\enspace
The Hermitian form $(\ , \ )_N$ is invariant under
 $\omega_{\frac{n-2}2, \epsilon, \Bbb R^n}(G)$.

Let $\Cal H$ be the completion of the pre-Hilbert space
 $(S (C_0^\infty(\Bbb R^n)), (\ , \ )_N)$.
 \newline
{\rm 3)}\enspace
 The Fr\'echet representation
 $\omega_{\frac{n-2}2, \epsilon, \Bbb R^n}$ of $G$
 on $\tilPsi(C^\infty(M)_\epsilon)$
 induces a representation of $G$ on the Hilbert space $\Cal H$, 
  which we shall denote by $(\minflat, \Cal H)$.
\newline
{\rm 4)}\enspace
The unitary representation $(\minflat, \Cal H)$ is irreducible.
\newline
{\rm 5)}\enspace
The map (see (2.8.7) for the definition)
$$
\alignat2
   \invtilPsi: & S (C_0^\infty(\Bbb R^n)) & \to & \Vpq
\\
               & \quad \cap                     & & \quad \cap
\\	   
               & C^\infty(\Bbb R^n)       & & C^\infty (M)    
\endalignat	       
$$
 extends uniquely to a unitary isomorphism
 between 
 $(\minflat, \Cal H)$ 
 and the minimal representation $(\varpi^{p,q}, \overline{\Vpq})$
 up to a scalar constant.
More precisely,
 for any $f_1, f_2 \in  \Cal H$,
 we have
$$
 (f_1, f_2)_N
 = 2^{2-n} (F_1, F_2)_M,
\tag \num.3
$$
 where  we put
$$
  F_i := \invtilPsi f_i,
 \quad (i=1, 2).
$$
\endproclaim
\remark{Remark}
1)\enspace
 We shall give a different proof of the first statement of Theorem~\ch.7
  in Theorem~\ch.9
 by using the Fourier transform of the Green kernel.
\newline
2)\enspace
 As we shall see in  the proof of Theorem~\num,
 $S(C^\infty_0(\Bbb R^n))$ is dense in $\tilPsi (\Vpq)$
 with respect to the above inner product.
\newline
3)\enspace
We can  realize naturally the Hilbert space $\Cal H$ as a subspace of
 the Schwartz distributions $\Cal S'(\Bbb R^n)$,
 namely,
 $\Cal F^{-1} T(L^2(C))$,
  as we shall prove in Theorem~\ch.9.
\endremark
\def\sc{8}
\sec{}
In order to prove Theorem~\ch.7, we need:
\proclaim{Lemma~\num}
Retain the notation of Theorem~\ch.7.
\newline
{\rm 1)}\enspace
For any $f_1, f_2 \in S(C^\infty_0(\Bbb R^n))$,
 we have
$$
 (f_1, f_2)_N = c_3
  ( F_1, F_2)_A,
\tag \num.1
$$
where 
$$
c_3:=
   \frac{\pi^{\frac{n}2} \ \Gamma(\frac{-1+n+\epsilon}{4})}
{\Gamma(\frac{n}2-1)}.
$$

We extend the Hermitian form $( \ , \ )_N$
from $S(C^\infty_0(\Bbb R^n))$ to 
$A^{\Bbb R^n}_{1, \epsilon} (\Psi_{\frac{n+2}2}^*(C^\infty(M)_\epsilon))$
by using the right side of (\num.1)
 (recall the inclusive relation (\ch.7.2)).
\newline
{\rm 2)}\enspace
The Hermitian form
$( \ , \ )_N$ is positive definite on
 $A^{\Bbb R^n}_{1, \epsilon} (\Psi_{\frac{n+2}2}^*(C^\infty(M)_\epsilon))$,
  in which
 $S (C_0^\infty(\Bbb R^n))$ is dense.
 In particular,
 the Hilbert space
 $\Cal H$ (see Theorem~4.7) coincides with the completion of the 
pre-Hilbert space 
$A_{1, \epsilon} (\Psi_{\frac{n+2}2}^*(C^\infty(M)_\epsilon))$.
\endproclaim

We first finish the proof of Theorem~\ch.7, and then
 give a proof of Lemma~\num.
\demo{Proof of Theorem~\ch.7}
(1) is clear from Lemma~\num~(1) and Lemma~\ch.4.
Next, let us prove (\ch.7.3).
We consider the $G$-intertwining operator
$$
\invtilPsi : \tilPsi \left(A^M_{1, \epsilon} (C^\infty(M)_\epsilon)\right)
 \to A_{1, \epsilon}^M(C^\infty(M)_\epsilon) \subset \Vpq,
$$
 or equivalently by (\ch.5.2),
$$
\invtilPsi :  A^{\Bbb R^n}_{1, \epsilon} 
\left(\Psi_{\frac{n+2}2}^* (C^\infty(M)_\epsilon)\right) \to  \Vpq.
$$

Combining (\ch.4.4) and (\num.1),
 we have
$$
(f_1, f_2)_N
 = c_3 (F_1, F_2)_A
 = c_1 c_3 (F_1, F_2)_M 
 = 2^{2-n} (F_1, F_2)_M,
$$
 where the second equality follows from
 a classical formula of the Gamma function:
$$
  2^{2 z-1} \Gamma(z) \ \Gamma(z+\frac12) = \sqrt{\pi} \ \Gamma(2 z).
$$
Thus, we have proved (\ch.7.3), especially,
$\invtilPsi$ is an isometry up to scalar.
Then $\invtilPsi$ extends to an isometric (up to scalar)
 $G$-intertwining operator
 $\Cal H \to \overline{\Vpq}$.
This map is surjective because $(\varpi^{p,q}, \overline{\Vpq})$
 is an irreducible unitary representation of $G$ (see Lemma~2.6).
All other statements  are now clear.
\qed
\enddemo
\demo{Proof of Lemma~\num}
1)\enspace
Suppose $f_i = S\varphi_i$, $\varphi_i \in C^\infty_0(\Bbb R^n)$ ($i = 1, 2)$.
Recall the constant $c_2$ in Proposition~\ch.6,
 we have
$$
\alignat1
  (f_1, f_2)_N
&= \int_{\Bbb R^n} (S \varphi_1)(x) \overline{\varphi_2(x)} d x
\\
& =
  \frac{1}{2 c_2}
  \int_{\Bbb R^n}
   (A^{\Bbb R^n}_{1, \epsilon} \varphi_1)(x) \ \overline{ \varphi_2(x)} \ d x
\\
& =
  \frac{1}{4 c_2}
  \int_{M}
   (A^M_{1, \epsilon}
 ({(\Psi^*_{\frac{n+2}2, \epsilon})^{-1}} \varphi_1)(b) \
 \overline{ {(\Psi^*_{\frac{n+2}2, \epsilon})^{-1}} \varphi_2 (b)} \ d b.
\tag  \num.2
\endalignat
$$
We put
$$
 \psi_i:=\frac{1}{c_2} 
 {(\Psi^*_{\frac{n+2}2, \epsilon})^{-1}} \varphi_i
\ \in C^\infty(M)_\epsilon, \quad
 (i=1,2).
$$
Then $F_i %
= A^M_{1, \epsilon} \psi_i$ by (\ch.5.2).
Therefore, 
$$
(\num.2) =
  \frac{1}{4 c_2} c_2^2
  \int_{M}
   (A^M_{1, \epsilon} \psi_1)(b) \
 \overline{ \psi_2 (b)} \ d b
 = \frac{\delta c_2}{4} (F_1, F_2)_A
$$
which equals the right side of (\num.1).
\newline
2) \enspace
It is enough to show that $S(C^\infty_0(\Bbb R^n))$ is dense
 in a pre-Hilbert space 
\newline
 $A_{1, \epsilon} (\Psi_{\frac{n+2}2}^*(C^\infty(M)_\epsilon))$
 because the inner product $(\  ,\ )_A$ is positive definite
 from Lemma~\ch.4.

Suppose
 $f = A_{1, \epsilon} \varphi$
  ($\varphi \in \Psi_{\frac{n+2}2}^*(C^\infty(M)_\epsilon)$)
 is orthogonal to $S(C^\infty_0(\Bbb R^n))$
 with respect to the inner product $(\  ,\ )_A$.
This means that $(f, A_{1, \epsilon} \phi)_A =0$
 for any $\phi \in C^\infty_0(\Bbb R^n)$.
Then we have
$$
 \int_M \invtilPsi f (b) \ \overline{\invtilPsi \phi(b)} \ d b = 0
\quad\text{ for any $\phi \in C^\infty_0(\Bbb R^n)$}.
$$
Since
 $\invtilPsi(C^\infty_0(\Bbb R^n))$ is dense in
 $L^2(M)_\epsilon := \set{f \in L^2(M)}{f(-u) = \epsilon f(u)}$,
 we have
 $\invtilPsi f  = 0$, and thus $f = 0$.
\qed
\enddemo
\def\sc{9}
\sec{}
We recall $\tilPsi(\Vpq) \subset \tilPsi(C^\infty(M)_\epsilon)
 \subset \Cal S'(\Bbb R^n)$,
 on which we can define the Fourier transform $\Cal F$.
We consider the following maps (see \S 4.3, \S 3.2, \S 3.4):
$$
 C_0^\infty(\Bbb R^n) 
 \overset{S}\to\rightarrow
 \tilPsi(\Vpq)
 \overset{\Cal F}\to\rightarrow
 \Cal S'(\Bbb R^n)
 \overset{T}\to\hookleftarrow
 L^2(C)
$$
Here is a description of the Fourier transform of the
 minimal unitary representation $(\minflat, \Cal H)$
 of $G = O(p,q)$
 which is obtained as the completion 
 of a pre-Hilbert space
  $(S (C^\infty_0(\Bbb R^{n})), (\ , \ )_N)$ (see Theorem~\ch.7).
\proclaim{Theorem~\num}
{\rm 1)}\enspace
$\Cal F S (C^\infty_0(\Bbb R^n))$ is contained in $T (L^2(C))$.
Furthermore,
it is dense in the Hilbert space $T (L^2(C))$.
\newline
{\rm 2)}\enspace
$T^{-1}\circ \Cal F \: S(C^\infty_0(\Bbb R^n)) \to L^2(C)$
 extends uniquely to a unitary map $\Cal H \to L^2(C)$
 up to a scalar constant.
 This constant is given explicitly by
$$
     (2 \pi)^n (f, f)_N = \| T^{-1} \circ \Cal F f\|_{L^2(C)}^2.
\tag \num.1
$$
\newline
{\rm 3)}\enspace
The $\overline{\Pmax}$-module $(\pi, L^2(C))$
 extends to an irreducible unitary representation of $G$, 
 denoted by the same letter $\pi$,
 so that $(2 \pi)^{-\frac{n}2}T^{-1}\circ \Cal F$ gives a unitary equivalence 
 between $(\minflat, \Cal H)$ and $(\pi, L^2(C))$.
\newline
{\rm 4)}\enspace
$T^{-1} \circ \Cal F \circ \tilPsi \: \Vpq \to L^2(C)$
 induces a unitary equivalence 
 between irreducible unitary representations
 $(\varpi^{p,q}, \overline{\Vpq})$ and $(\pi, L^2(C))$,
 up to a scalar constant given by:
$$
 (\phi, \phi)_M = 
 4 \pi^n
\| T^{-1}\circ \Cal F \circ \tilPsi \phi\|^2_{L^2(C)}
\tag \num.2
$$
\endproclaim
\demo{Proof}
If $f = E_0 * \varphi$ ($\varphi \in C_0^\infty (\Bbb R^n)$), 
 then it follows from the integration formula of the Green function
 (see Proposition~4.2) that its Fourier transform is given by 
$$
\Cal F f = \Cal F(E_0 * \varphi) = (\Cal F E_0) (\Cal F \varphi)
   = (\Cal F \varphi) \delta (Q) = T((\Cal F \varphi)|_C).
$$
Since $\varphi \in C_0^\infty(\Bbb R^n)$,
 we have $(\Cal F \varphi)|_C \in L^2(C)$. 
Hence,
 $\Cal F S (C^\infty_0(\Bbb R^n))$ is contained in $T(L^2(C))$.
Then we have (\num.1), 
 as follows from
$$
\| (\Cal F \varphi)|_C \|_{L^2(C)}^2
=   ((\Cal F E_0) (\Cal F \varphi), \Cal F \varphi) 
=   (\Cal F (E_0 * \varphi), \Cal F \varphi) 
= (2 \pi)^n (f, f)_N.
\tag \num.3
$$
We note that (\num.3) gives
 a different proof that $(\ , \ )_N$ is a positive definite Hermitian
 form on $S(C_0^\infty(\Bbb R^n))$ (see Remark after Theorem~\ch.7).

It follows from Lemma~3.4 that  $T^{-1} \circ \Cal F$
 is an $\overline\Pmax$-intertwining operator from
 $(\varpi_{\frac{n-2}2, \epsilon, \Bbb R^{n}}|_{\overline\Pmax},
 S(C^\infty_0(\Bbb R^n)))$
 to $(\pi, L^2(C))$.
This map is isometric up to a scalar by (\num.1).
Then, it extends naturally to an 
 $\overline\Pmax$-intertwining operator
 from $(\minflat|_{\overline\Pmax}, \Cal H)$ to $(\pi, L^2(C))$,
 which is surjective because
 $(\pi, L^2(C))$ is irreducible (see Proposition~3.3).
Hence we have proved (1), (2) and (3).
The statement (4) follows from (2) and Theorem~\ch.7.
Thus, we have proved Theorem~\ch.9.
\qed
\enddemo
\def\sc{10}
\sec{}
By semisimple theory,
 it is known 
 that a minimal representation is still irreducible 
 when restricted to any maximal parabolic subgroup.
 In particular,
  $\varpi^{p,q}$ is irreducible as an $\overline{\Pmax}$-module.
In our case, this fact can be strengthened as follows:
\proclaim{Corollary~\num}
The restriction of the minimal representation $\varpi^{p,q}$ 
 (equivalently, $\minflat$) 
to %
 $\Mmax_+ \overline{\Nmax} \simeq O(p-1, q-1) \times \Bbb R^{p+q-2}$
 is still irreducible.
\endproclaim
\demo{Proof}
Theorem~\ch.7 and Theorem~\ch.9 show that
$\varpi^{p,q}$, $\minflat$ and $\pi$ are unitary equivalent to one another.
Now, Corollary follows from Proposition~3.3. 
\qed
\enddemo
\def\sc{11}
\sec{}
For the convenience of the reader, we summarize the maps used in the proofs.
$$
\alignat4
& C^\infty(M) & \overset{\tilPsi}\to\longrightarrow
             & C^\infty(\Bbb R^n) \cap \Cal S'(\Bbb R^n)
             & \overset{\Cal F}\to\longrightarrow
             & \ \Cal S'(\Bbb R^n) && %
\tag \num.1
\\
&\quad \cup  && \qquad\quad \cup &&\quad \cup &&
\\
              &\quad \Vpq
              &\ \rarrowsim \ 
              &\tilPsi (\Vpq)
              &\ \rarrowsim \
              &\Cal F \tilPsi (\Vpq)
               &&
\\
              &\quad \cap\text{dense}
              & 
              &\quad \cap\text{dense}
              &
              &\quad \cap\text{dense}
	      &
\\
              &(\varpi^{p,q}, \overline{\Vpq})
              &\ \rarrowsim \ 
              &\quad (\minflat, \Cal H)
              &\ \rarrowsim \
              &T(L^2(C))
              &\underset{T}\to{\larrowsim} 
	      & (\pi, L^2(C)).
\endalignat
$$             
In the last line, 
 we have written also the notation for unitary representations.

\def\ch{5}
\head
\S \ch.  
   Bessel function and an integral formula of spherical functions 
\endhead
\def\sc{1}
\sec{}
In this section, we shall compute explicitly the lowest $K$-type
of our minimal representation in the $N$-picture, i.e. find it
as a solution to $\fLap f = 0$, and also its Fourier transform
as a function on $C$; this turns out to be written
 in terms of a Bessel function (see Theorem~\ch.5). 
Note that except when $p = q$ we are not dealing with a spherical
representation of $G$ (namely, there is no non-zero $K$-fixed vector
 in our representation). 
At the end of this section, we reformulate
the equivalent realizations of the minimal representation found
in the previous section, using now the minimal $K$-type to
understand the different pictures. 
\def\sc{2}
\sec{}
Without loss of generality,
 we may and do assume $p \ge q$ in this section.
Instead of $K$-fixed vectors,
 our idea here is to focus on an $O(p) \times O(q-1)$-fixed vector.
Then,
 it follows from Lemma~2.6 that such a vector, 
 which we shall denote by $F_0$,
 is unique in our minimal representation $(\varpi^{p,q}, \Vpq)$ up to 
 a scalar multiple,
 and is contained in the minimal $K$-type
 of the form $\boldkey 1 \boxtimes \spr{\frac{p-q}2}{q}$.
(We note that this $K$-type is not one dimensional if $p \neq q$).

We shall find an explicit formula of the Fourier transform
 of this vector $F_0$ after a conformal change of coordinates.
We start with the following classical lemma,
 for which we give a proof for the sake of the completeness.
We take a coordinate $(u_{p+1}, \dots, u_{p+q})$ in $\Bbb R^q$
 and realize $O(q-1)$ in $O(q)$ such that
 it stabilizes the last coordinate $u_{p+q}$.
\proclaim{Lemma~\num}
For any $l \in \Bbb N$,
 $O(q-1)$-invariant spherical harmonics of degree $l$ form
 a one dimensional vector space.
More precisely,
 we have
$$
   \spr{l}{q}^{O(q-1)}
  \simeq
   \Bbb C \ \trF(\frac{-l}2, \frac{q-2+l}2; \frac{q-1}2;
             {u_{p+1}}^2 + \dots + {u_{p+q-1}}^2).
$$
\endproclaim
\demo{Proof}
In terms of the polar coordinate of $S^{q-1}$:
$$
\Phi_{++} : S^{q-2} \times (0, \frac{\pi}2) \to S^{q-1},
 \ 
 (y, \theta) \mapsto ((\sin \theta)y, \cos \theta),
\tag \num.1
$$
 the Laplace-Beltrami operator on $S^{q-1}$ takes the form:
$$
   \Delta_{S^{q-1}} = \frac{\partial^2}{\partial \theta^2}
                     + (q-2) \cot \theta \frac{\partial}{\partial \theta}
              +\frac{1}{\sin^2 \theta} \Delta_{S^{q-2}}.
\tag \num.2
$$
If $F \in \spr{l}{q}$ is $O(q-1)$-invariant,
 then 
 $F \circ \Phi_{++}(\theta, y)$ depends only on $\theta$,
 for which we write $h(\theta)$.
Then $h(\theta)$ is an even function satisfying:
$$
 (\frac{d^2}{d \theta^2}
     + (q-2) \cot \theta \frac{d}{d \theta}
     + l(l+q-2)) h(\theta) = 0.
$$
Since $h(\theta)$ is regular at $\theta = 0$,
 it is a scalar multiple of the Jacobi function:
$$
\alignat1
  \varphi^{\frac{q-3}2, -\frac12}_{\sqrt{-1}(l+\frac{q-2}2)}(\sqrt{-1} \theta)
&=
   \trF(\frac{-l}2, \frac{q-2+l}2; \frac{q-1}2; \sin^2 \theta)
\\
&=
   \trF(-l, q-2+l; \frac{q-1}2; \frac{1-\cos\theta}{2}).
\endalignat
$$
Thus, we have proved the lemma.
\qed
\enddemo

\def\sc{3}
\sec{}
In view of Lemma~2.6, 
 the special case of Lemma~\ch.2 with $l= \frac{p-q}2$ yields:
\proclaim{Proposition~\num}
Suppose $p \ge q \ge 2$, $p+q \in 2 \Bbb N$ and $(p,q) \neq (2,2)$.
Let $(u_1, \dots, u_{p+q})$ be the coordinate of $M = S^{p-1} \times S^{q-1}$
 in $\Bbb R^{p+q}$.
We define a function $F_0 : M \to \Bbb C$ by
$$
F_0(u_1, \dots, u_{p+q}) 
:= \trF(\frac{q-p}4, \frac{p+q-4}4; \frac{q-1}2;
 {u_{p+1}}^2+\dots+{u_{p+q-1}}^2).
\tag \num.1
$$
Then $F_0$ is an $O(p) \times O(q-1)$-invariant analytic function on $M$
 satisfying the Yamabe equation
$
   \tilLap{M} F_0 = 0.
$
Conversely, any such function is a scalar multiple of $F_0$.
\endproclaim
\remark{Remark}
If $p+q$ is odd, then  $F_0$ in the right side of (\num.1)
 still  gives a solution $\tilLap{M} F = 0$
 on an open dense set of $M$ such that $u_{p+q}\neq 0$.
Furthermore, $F_0$ is a continuous function on $M$.
However, it does not solve the Yamabe equation as a distribution on $M$. 
\endremark

\def\sc{4}
\sec{}
We recall $\tau(z', z'')$ is a conformal factor defined in (2.8.1).
Let us define $f_0 := \tilPsi F_0$, namely,
$$
f_0(z',z'') :=    \tau(z', z'')^{-\frac{p+q-4}2}
  \trF (\frac{q-p}4, \frac{p+q-4}4; \frac{q-1}2;
 \frac{|z''|^2}{\tau(z', z'')^2}).
\tag \num.1
$$
We note that $|\tau(z', z'')| \ge | z''|$ 
 for any $(z', z'') \in \Bbb R^{p-1,q-1}$.
The equality holds if and only if 
$|z'|^2 -|z''|^2 = -4$.

The following Proposition is immediate from Lemma~2.8 and Proposition~\ch.3:
\proclaim{Proposition~\num}
With the same assumption on $p, q$ in Proposition~\ch.3, we have:
\newline
\rm{1)}\enspace
$f_0$ is a real analytic function on $\Bbb R^n$ that
 solves $\fLap f_0 = 0$.
\newline
\rm{2)}\enspace
$f_0$ is $O(p) \times O(q-1)$-invariant.
\endproclaim

We say $F_0$ is the {\it generating function} of $\Vpq = \Ker \tilLap{M}$,
 and $f_0$ is that of $(\minflat, \Cal H)$.

\remark{Remark}
More strongly than Proposition~\num,
 one can prove that
 $f_0$ is a real analytic solution of $\fLap f_0 = 0$ if $p+q > 4$
 by using Proposition~5.6, where we do not assume that
 $p+q$ is even.
The real analyticity is not obvious from the expression (\num.1)
 in the neighbourhood of the hypersurface of
 $|z'|^2 -|z''|^2 = -4$.
\endremark

\def\sc{5}
\sec{}
We recall the definitions of Bessel functions:
$$
\alignat2
  J_\nu(z) &= \sum_{m=0}^\infty 
\frac{(-1)^m (\frac12 z)^{\nu+2m}}{m!\ \Gamma(\nu+m+1)}
\qquad
&&\text{(Bessel function)},
\\
  I_\nu(z)
 &= \sum_{m=0}^\infty \frac{(\frac12 z)^{\nu+2m}}{m!\ \Gamma(\nu+m+1)}
\qquad
&&\text{(modified Bessel function of the first kind)},
\\
  K_\nu(z) &= \frac{\pi}2 \frac{I_{-\nu}(z)-I_\nu(z)}{\sin \nu\pi}
\qquad
&&\text{(modified Bessel function of the second kind)}.
\endalignat
$$
Then $K_\nu$ satisfies
$$
   \left( z^2\frac{d^2}{d z^2} + z \frac{d}{d z} - (z^2 + \nu^2)\right)
   K_\nu(z) = 0.
$$

The asymptotic of the functions $K_\nu$  is well-known, for example
$z^{-\nu} K_\nu (z)$ decays exponentially as $z \to +\infty$.

\proclaim{Theorem~\num}
We put $|\zeta| := ({\zeta_1}^2 + \dots + {\zeta_{n}}^2)^{\frac12}$
 for $\zeta = (\zeta_1, \dots, \zeta_{n}) \in \Bbb R^n$.
Let $F_0 \in C^\infty(M)$ be the generating function
 of $\Vpq = \Ker \tilLap{M}$ (see Proposition~\ch.3).
Then
$$
(\Cal F \tilPsi F_0)(\zeta)
=
(2\pi)^{\frac{p+q-2}{2}} 2^{-\frac{p-5}2}
 \frac{\Gamma(\frac{q-1}2)}{\Gamma(\frac{p+q-4}2)}
 |\zeta|^{\frac{3-q}2} K_{\frac{q-3}2}(2 |\zeta|) \delta(Q).
$$
\endproclaim

Note that the
 $|\zeta|^{\frac{3-q}2} K_{\frac{q-3}2}(2 |\zeta|)$
 belongs to the Hilbert space $L^2(C)$
 if $p+q > 4$ by
the asymptotic behaviour of the Bessel function $K_\nu$
 and by the explicit form of $\delta(Q)$ in (3.3.3).
 
\def\sc{6}
\sec{}
Theorem~\ch.5 follows from the following Proposition:
\proclaim{Proposition~\num}
We write $\Cal F^{-1}$ for the inverse Fourier transform. 
With notation in (\ch.4.1), we have
$$
 \Cal F^{-1} 
 \left(|\zeta|^{\frac{3-q}2} K_{\frac{q-3}2}(2 |\zeta|) \delta(Q)\right) (z)
=
 \frac{\Gamma(\frac{p+q-4}2)}{2^{\frac{q+3}2} \pi^{\frac{p+q-2}{2}}
 \Gamma(\frac{q-1}2)} f_0(z).
\tag \num.1
$$
\endproclaim
\demo{Proof}
Let $\phi(r)$ be a function of one variable,
 which will be taken later to be $r^{\frac{3-q}2} K_{\frac{q-3}2}(2 r)$.
Then, it follows from (3.3.3) that for $z = (z',z'') \in \Bbb R^{p-1,q-1}$,
$$
\alignat1
&2\Cal F^{-1} \left(\phi(|\zeta|) \delta(Q)\right)(z)
\\
&=
(2 \pi)^{-(p+q-2)} \int_0^\infty \int_{S^{p-2}} \int_{S^{q-2}}
 \phi(r) e^{-\sqrt{-1} ((z', r \omega)+(z'', r \eta))} \ r^{p+q-5}
  d r \ d \omega \ d \eta
\intertext{
   Using the formula
$
            \int_{S^{m-1}} e^{\sqrt{-1} t (\eta, \omega)} d \omega
            =
            (2 \pi)^{\frac{m}2} t^{1 - \frac{m}2} J_{\frac{m}2 -1}(t),
$
we have
}
=
&
 (2\pi)^{-\frac{p+q-2}{2}}
\int_0^\infty \phi(r) \
 (r |z'|)^{\frac{3-p}2} J_{\frac{p-3}2}(r |z'|)
 (r |z''|)^{\frac{3-q}2} J_{\frac{q-3}2}(r |z''|)
 r^{p+q-5} \ dr
\\
=
&
 (2\pi)^{-\frac{p+q-2}{2}}
  |z'|^{\frac{3-p}2} |z''|^{\frac{3-q}2} 
\int_0^\infty \phi(r)
  J_{\frac{p-3}2}(r |z'|) J_{\frac{q-3}2}(r |z''|)
 r^{\frac{p+q-4}2} \ dr.
\endalignat
$$
Now, put
 $\phi(r) := r^{\frac{3-q}2} K_{\frac{q-3}2}(2 r)$.
We use the following formula of the Hankel transform due to Bailey
 \xbail\ (see also \xerdIntII, \S 19.6 (8))
$$
  \align
     &\int_0^{\infty}
     t^{\lambda-1}J_{\mu}(a t) J_{\nu}(b t) K_{\rho}(c t) d t
\\
    =& \frac {2^{\lambda-2}
              a^{\mu} b^{\nu}
              \Gamma({\frac 1 2}(\lambda+ \mu + \nu-\rho))
              \Gamma({\frac 1 2}(\lambda+ \mu + \nu+\rho))}
              {c^{\lambda+ \mu + \nu} \Gamma (\mu+1)\Gamma(\nu+1)}
\\
     &\times
       F_4(\frac 1 2 (\lambda+ \mu + \nu -\rho),
           \frac 1 2 (\lambda+ \mu + \nu +\rho);
            \mu+1, \nu+1;
            -\frac {a^2}{c^2}, -\frac {b^2}{c^2}).  
  \endalign
$$
Here, $F_4$ is the Appell hypergeometric function of two variables,
 defined by
$$
F_4(a, b; c, d; x, y)
=
\sum_{i=0}^\infty \sum_{j=0}^\infty
\frac{(a)_{i+j} (b)_{i+j}}{i !  j ! (c)_i (d)_j} x^i y^j.
$$
Then we have:

$$
\alignat1
&\Cal F^{-1} \left(\phi(|\zeta|) \delta(Q)\right)(z)
\\
=
&
 (2\pi)^{-\frac{p+q-2}{2}}
  |z'|^{\frac{3-p}2} |z''|^{\frac{3-q}2} 
\int_0^\infty  K_{\frac{q-3}2}(2  r)
  J_{\frac{p-3}2}(r |z'|) J_{\frac{q-3}2}(r |z''|)
 r^{\frac{p-1}2} \ d r
\\
=
&
 (2\pi)^{-\frac{p+q-2}{2}} 2^{\frac{p-3}2}
 \frac{\Gamma(\frac{p+q-4}2)}{\Gamma(\frac{q-1}2)}
F_4(\frac{p-1}2, \frac{p+q-4}2; \frac{p-1}2, \frac{q-1}2; 
\frac{-|z'|^2}{4};\frac{-|z''|^2}{4}).
\endalignat
$$
\enddemo
\def\sc{7}
\sec{}
Then,
 the proof of Proposition~\ch.6
 will be finished by showing the following reduction formula:
\proclaim{Lemma~\num}
Let $\tau(z', z'')$ be the conformal factor defined in (2.8.1).
We have
$$
\alignat1
&F_4(\frac{p-1}2, \frac{p+q-4}2; \frac{p-1}2, \frac{q-1}2; 
\frac{-|z'|^2}{4};\frac{-|z''|^2}{4})
\\
=
&
 \tau(z', z'')^{-\frac{p+q-4}2}
\trF
(\frac{q-p}4, \frac{p+q-4}4; \frac{q-1}2; \frac{|z''|^2}{\tau(z', z'')^2}).
\tag \num.1
\endalignat
$$
\endproclaim
\demo{Proof of Lemma~\num}
We recall a reduction formula of Appell's hypergeometric functions 
(see \xerdHigI, \S 5.10, (8)):
$$
\alignat1
&F_4(\alpha, \beta; 1+ \alpha-\beta, \beta;
 \frac{-x}{(1-x)(1-y)},\frac{-y}{(1-x)(1-y)})
\tag \num.2
\\
=
&(1-y)^\alpha \trF(\alpha, \beta; 1+\alpha-\beta; \frac{-x(1-y)}{1-x})
\endalignat
$$
and a quadratic transformation for hypergeometric functions
(see \xerdHigI, \S 2.11 (32)):
$$
\trF(\alpha, \beta, 1+\alpha-\beta; z)
= (1-z)^{-\alpha} 
\trF(\frac{\alpha}2, \frac{\alpha + 1 - 2 \beta}2; 1+\alpha -\beta; \frac{-4 z}{(1-z)^2}).
\tag \num.3
$$
Combining (\num.2) with (\num.3) for
 $\alpha = \frac{p+q-4}2$ and $\beta = \frac{p-1}2$,
 and using the symmetry of $a$ and $b$; $(c,x)$ and $(d, y)$
 in $F_4(a,b;c,d; x,y)$,
 we have
$$
\alignat1
&F_4(\frac{p-1}2, \frac{p+q-4}2; \frac{p-1}2; \frac{q-1}2;
 \frac{-x}{(1-x)(1-y)},\frac{-y}{(1-x)(1-y)})
\\
&=
 \left(\frac{(1-x)(1-y)}{1-x y}\right)^{\frac{p+q-4}2}
 \trF(\frac{p+q-4}4, \frac{q-p}4; \frac{q-1}2; \frac{4 y (1-x)(1-y)}{(1- x y)^2})
\endalignat
$$

If we put
$$
\frac{|z'|^2}4 = \frac{x}{(1-x)(1-y)},
\quad
\frac{|z''|^2}4 = \frac{y}{(1-x)(1-y)}
$$
then a simple computation shows
$$
 \tau(z', z'')^2 =  \left(\frac{1-x y}{(1-x)(1-y)}\right)^2,
\quad
\frac{|z''|^2}{\tau(z',z'')^2} =  \frac{4 y (1-x)(1-y)}{(1- x y)^2}.
$$
Thus, Lemma~\num\ is proved.
\qed
\enddemo

\def\sc{8}
\sec{}
Using the generating function $F_0$ we may recover the whole representation
by letting the Lie algebra of $G$ act. Let us see how our previous
results may be reformulated:  
It follows from the definition of $\widehat{\varpi}_{\frac{n-2}2, \epsilon}$
 (see \S 3.2) that the linear map
$$
   \Cal F \circ \tilPsi \: C^\infty(M) \to {\Cal S}'(\Bbb R^n)
$$
 induces a natural intertwining map from $(\varpi^{p,q}, \Vpq)$
 to
 $(\widehat\varpi_{\frac{n-2}2, \epsilon}, \Cal F \tilPsi (\Vpq))$
 as $G$-modules and also as $\frak g$-modules.
Here again $\epsilon = (-1)^{\frac{p-q}2}$,
 $\Vpq = \Ker\tilLap{M} \ (\subset C^\infty(M))$ and 
 $\tilPsi (\Vpq) \subset \Ker \fLap$.

For $b \in \Bbb R^n, m \in O(p-1,q-1)$,
we put a function on the cone $C$ by
$$
    \psi_{b, m} (\zeta) :=
e^{\sqrt{-1}\langle b, \zeta \rangle}
 |m \zeta|^{\frac{3-q}2} K_{\frac{q-3}2}(2 |m \zeta|).
\tag \num.1
$$
In particular,
 we have
$$
   \psi_{0, e}(\zeta) = |\zeta|^{\frac{3-q}2} K_{\frac{q-3}2}(2 |\zeta|).
$$ 

Here we give explicit functions which are dense in the minimal 
 representations.
\proclaim{Theorem~\num}
Suppose $p \ge q \ge 2$, $p + q \in 2 \Bbb N$, and $(p,q) \neq (2,2)$.
\newline
{\rm 1)}\enspace
$\psi_{0, e}(\zeta)$ is a $K$-finite vector of $(\pi, L^2(C))$.
It belongs to the minimal $K$-type of $\pi$.
\newline
{\rm 2)}\enspace
 $\Bbb C\text{-span}\set{\psi_{b,m}}{
 b \in \Bbb R^n, m \in O(p-1,q-1)}$
 is a dense subspace of the minimal representation
$(\pi, L^2(C))$.
\newline
{\rm 2$'$)}\enspace
 $\Bbb C\text{-span}\set{\psi_{b,m}\delta(Q)}{
 b \in \Bbb R^n, m \in O(p-1,q-1)}$
 is a dense subspace of the minimal representation
$(\widehat\varpi_{\frac{n-2}2, \epsilon}, \Cal F \tilPsi(\overline\Vpq))$.
\endproclaim
\demo{Proof}
This follows by combining Theorem 4.9 and Lemma 5.5, and using that
by Mackey theory we have an irreducible representation of the parabolic
group.
\qed
\enddemo 

\def\sc{9}
\sec{}
The advantage in the realization on $L^2(C)$ is
 that the action for $\overline\Pmax$ and the inner product are
 easily described.
On the other hand,
 the action of $K$ is not easy to be described,
 and especially, the $K$-finiteness in the statement (2) is non-trivial.

Let $U(\frak g)$ the enveloping algebra of
 the complexified Lie algebra $\frak g$.
We define a subspace of $\Cal S'(\Bbb R^n)$ by
$$
 U := d \widehat\varpi_{\frac{n-2}2, \epsilon} (U(\frak g))
 (\psi_{0, e} \delta(Q)).
\tag \num.1
$$
That is, $U$ 
 is the linear span of a Bessel function
 $|\zeta|^{\frac{3-q}2} K_{\frac{q-3}2}(2 |\zeta|) d \mu$ 
 on the cone $C$ and 
 its iterative differentials corresponding to 
 the action of the Lie algebra $\frak g$.

We have seen that $\Cal F \tilPsi(\Vpq) \subset T(L^2(C))$
 in Theorem 4.9.
We may restate in this way:
\proclaim{Theorem~\num}\enspace
{\rm 1)}\enspace
$U$ is an infinitesimally unitary
 $(\frak g, K)$-module via $\widehat\varpi_{\frac{n-2}2, \epsilon}$.
\newline  
{\rm 2)}\enspace
 $U$ is dense in the Hilbert space $T(L^2(C))$.
\newline  
{\rm 3)}\enspace
The completion of (1) defines an irreducible unitary representation
 of $G$ on $T(L^2(C))$, and then also on $L^2(C)$.
This gives an extension of $\pi$ from $\overline{\Pmax}$ to $G$.
\endproclaim
This has already been done by Theorem~4.9 and the irreducibility
 of the minimal representation.

One of non-trivial parts of the above assertion is to show
$$
   U \cap T(L^2(C)) \neq \{0\}
$$
which was proved in Theorem~5.5.

\def\ch{6}
\head
\S \ch. 
 Explicit inner product on solutions $\fLap f =0$
\endhead

\redefine\sc{1}
\sec{}
The aim of this section is to provide an explicit inner product
 on a certain subspace (see (\ch.2.1))
 of solutions of the ultrahyperbolic equation $\fLap f = 0$,
 such that its Hilbert completion gives 
 the unitarization of the  minimal representation of $O(p,q)$.

Roughly speaking, 
 the inner product will be given in terms of the integration over
 a hyperplane after convoluting a distribution along the normal direction.

We assume $n > 2$.
We fix $i \in \{1, 2, \dots, n\}$ once and for all.
The hyperplane on which we integrate will be
 $\set{z =(z_1, \dots, z_n) \in \Bbb R^n}{z_i= 0}$,
 for which we simply write $\Bbb R^{n-1}$.

\redefine\sc{2}
\sec{}
Let $C^{(i)} := \set{\zeta \in C}{\zeta_i \neq 0}$,
 an open dense subset of the null cone $C$ (see (3.3.2)).
We note that $C_0^\infty(C \setminus \{0\})$ is dense in $L^2(C)$.
We define a subspace of solutions of $\Ker \fLap$ by
$$
    W := \Cal F^{-1} \circ T (C_0^\infty(C^{(i)})).
\tag \num.1
$$
Here,
 we recall $T : L^2(C) \hookrightarrow \Cal S'(\Bbb R^n)$
 is the embedding via the measure $d \mu$ on the cone $C$.
By the Paley-Wiener theorem for compactly supported distributions,
 $W$ consists of real analytic solutions of $\fLap f = 0$.

Using an interpretation of the Dirac delta function in terms of hyperfunctions:
$$
\delta(z_i)
 = \frac{1}{2 \pi \sqrt{-1}}
     \left(
            \frac{1}{z_i - \sqrt{-1} 0} - \frac{1}{z_i + \sqrt{-1} 0}
     \right),
$$
 we decompose $f \in W$ as
$$
     f(z) = f_+^{(i)}(z) + f_-^{(i)}(z),
\tag \num.2
$$ 
 where $\pmi{f}(z)$ is defined by the convolution in the $z_i$-variable:
$$
  \pmi{f}(z) := \frac{1}{2 \pi \sqrt{-1}} \cdot
                \frac{\pm 1}{z_i \mp \sqrt{-1} 0}
                     * f(z).
\tag \num.3
$$
We shall see later that the decomposition (\num.2) makes sense
 not only for $f \in W$ but also
 for any $f \in \Cal F^{-1} \circ T (L^2(C))$ (see Lemma~\ch.5).
We set 
$$
   (f, f)_W := 
\frac{1}{\sqrt{-1}}
\int_{\Bbb R^{n-1}}
 \left( 
  f_+^{(i)} \overline{\frac{\partial f_+^{(i)}}{\partial z_i }} 
  - 
  f_-^{(i)} \overline{\frac{\partial f_-^{(i)}}{\partial z_i }}
 \right)|_{z_i = 0}
\,
  d z_1 \dotsb \widehat{d z_i} \dotsb d z_n
\tag \num.4
$$
\proclaim{Theorem~\num}
Fix  any $i \in \{1,2,\dots, n\}$.
\newline
{\rm 1)}\enspace
The formula (\num.4) defines a positive definite Hermitian form on $W$,
 a subspace of solutions of the ultrahyperbolic operator
  $\fLap = \frac{\partial^2}{{\partial z_1}^2} + \dots +
    \frac{\partial^2}{{\partial z_{p-1}}^2}
    - \frac{\partial^2}{{\partial z_p}^2} - \dots -
    \frac{\partial^2}{{\partial z_{p+q-2}}^2}$ on $\Bbb R^{p-1,q-1}$.
\newline
{\rm 2)}\enspace
The inner product (\num.4) 
 is independent of 
 the choice of $i \in \{1,2,\dots,n\}$.
\newline
{\rm 3)}\enspace
The action of $G = O(p,q)$ preserves the inner product (\num.4),
 so that
 the Hilbert completion $\overline{W}$ of $W$ 
 defines a unitary representation of
 $G=O(p,q)$.
\newline
{\rm 4)}\enspace
The resulting unitary representation
 is unitarily equivalent to the minimal representation
 $(\varpi^{p,q}, \overline{V^{p,q}})$.
The $G$-intertwining operator 
 $\tilPsi \: \overline{V^{p,q}} \to \overline{W}$
  gives a unitary equivalence up to a scalar constant.
\newline
{\rm 5)}\enspace
$c \Cal F^{-1} \circ T : L^2(C) \to \overline{W}$
 is a unitary $G$-intertwining operator,
 if we put $c = 2^{\frac{n+2}2}\pi^{\frac{n+1}{2}}$.
\endproclaim

This Theorem gives a new formulation of the Hilbert space of the
minimal representation purely in terms of intrinsic objects in
the flat space $\Bbb R^{p-1,q-1}$ where the differential equation
is the classical ultrahyperbolic one. It generalizes the $q = 2$
case where an inner product was known in terms of integration
of Cauchy data - here one could interpret the inner product in terms
of the energy generator. The interesting property about the inner product
is its large invariance group; even translational invariance amounts
to a remarkable \lq\lq conservation law", and we may also note that the
integration over a coordinate hyperplane can be replaced by  
integration over any non-characteristic hyperplane
(since such a hyperplane is conjugate to
 either $z_1=0$ or $z_n=0$ by $O(p-1,q-1) \ltimes \Bbb R^{p-1,q-1}$), or even
 the image of such a hyperplane under conformal inversion.

The strategy of the proof of Theorem~\ch.2 is as follows:
We recall from Theorem~4.7 that the Hilbert space
 $\Cal H$ is the completion of the space $S(C^\infty_0(\Bbb R^n))$
  with respect to another inner product $(\ , \ )_N$
   ($S$ is an integral transform by the Green kernel).
Since $C^\infty_0(C \setminus \{0\})$ is dense in $L^2(C)$
 and since $T^{-1} \circ \Cal F : \Cal H \to L^2(C)$ is an isomorphism
 of Hilbert spaces (up to a scalar) by Theorem~4.9~(1),
  $W$ is a dense subspace of the Hilbert space $(\Cal H, (\ , \ )_N)$.
In light of this,
 the key ingredient of the proof of Theorem~6.2 is 
 to give a formula of $(\ , \ )_W$ by means of $(\ , \ )_{L^2(C)}$.
We shall prove:
 $$
    2 (2 \pi)^{n+1} (f, f)_W
     =
     \| T^{-1}\circ \Cal F f \|_{L^2(C)}^2
 \quad
 \text{for any } f \in W.
 \tag \num.5
 $$
Once we prove (\num.5), it follows from (4.9.1) that
 $$
   4 \pi (\ , \ )_W
    = (\ , \ )_N
 $$
 on the subspace $W$.
 In particular, we have
 $$
 \Cal H = \overline{W},
 $$
and all other statements of Theorem~\ch.2 on 
our inner product $( \ ,  \ )_W$ are clear from
 the corresponding results on the inner product $(\ , \ )_N$
  proved in Theorem~4.7 and Theorem~4.9
  (e.g. Theorem~6.2~(4) corresponds to Theorem~4.7~(5);
  Theorem~6.2~(5) to Theorem~4.9~(3)).
 
The rest of this section is devoted to the proof of the formula (\num.5).

\def\sc{3}
\sec{}
We define an open subset of the cone $C$ by
$$
       \pmi{C} := \set{\zeta \in C}{\pm \zeta_i > 0}.
$$
Then 
$$
    C^{(i)} := C_+^{(i)} \cup C_-^{(i)}
$$ 
is an open dense subset of the cone $C$,
 and we have a direct sum decomposition of the Hilbert space:
$$
    L^2(C) = L^2(C_+^{(i)}) \oplus L^2(C_-^{(i)}).
$$
We define
 the Heaviside function $\pmi{Y}(\zeta)$ of the variable $\zeta_i$ by
$$
   \pmi{Y}(\zeta_1, \dots, \zeta_n) 
       = \cases
             1 & \text{if $ \pm \zeta_i > 0$},
          \\
             0 & \text{if $ \pm \zeta_i \le 0$},
         \endcases
\qquad \text{ for $\zeta = (\zeta_1, \dots, \zeta_n) \in \Bbb R^n$}.
$$

For $\phi \in C_0^\infty(C^{(i)})$,
 we put 
$$
 \pmi{\phi} := \pmi{Y} \phi.
$$
Then $\operatorname{Supp}{\pmi{\phi}} \subset \pmi{C}$,
 and we have 
$$
\alignat1
\phi &= \phi_+^{(i)} + \phi_-^{(i)},
\tag \num.1
\\
\|\phi\|^2_{L^2(C)}
 &= \|\phi_+^{(i)}\|^2_{L^2(C)} + \|\phi_-^{(i)}\|^2_{L^2(C)}.
\tag \num.2
\endalignat
$$

\redefine\sc{4}
\sec{}
Let us take the $n-1$ variables
 $\zetaprime$
 as a coordinate on $\pmi{C}$.
Then we have
$$
  \zeta_i 
  = \pm 
  \sqrt{Q^{(i)}(\zetaprime)}
\tag \num.1
$$
 on $\pmi{C}$, respectively, 
 if we put
$$
  Q^{(i)}(\zetaprime)
 := -\varepsilon_i \left( {\zeta_1}^2 + {\zeta_2}^2 + \dotsb \widehat{\pm {\zeta_i}^2} \pm \dotsb
 - {\zeta_{n-1}}^2 - {\zeta_n}^2\right),
 \tag \num.2
$$
where $\varepsilon_i = \pm 1$ is the signature of ${\zeta_i}^2$ in the quadratic form $Q(\zeta)$
 as in (2.5.1).
We note that $Q^{(i)}(\zetaprime) \ge 0$ on the cone $C$,
 and the map
$$
  (\zetaprime) \mapsto (\zeta_1, \dots, \pm \sqrt{Q^{(i)}}, \dots, \zeta_n)
$$
gives a bijection from
$\set{(\zetaprime)\in\Bbb R^n}{Q^{(i)}(\zetaprime) \neq 0}$ onto $C^{(i)}$.
By substituting (\num.1) into $\pmi{\phi}$, we put
$$
  \pmi{\varphi}(\zetaprime)
  := 
    \pmi{\phi}(\zeta_1, \dots, \pm \sqrt{Q^{(i)}}, \dots, \zeta_n).
$$
Since the measure $d \mu$ on the cone $C$ is of the form
$
\frac{1}{2\sqrt{{Q^{(i)}}}}
 d\zeta_1 \dotsb \widehat{d \zeta_i} \dotsb d \zeta_n,
$
we have
$$
\alignat1
 \|\pmi{\phi}\|^2_{L^2(C)}
  &=\int_{\Bbb R^{n-1}} 
\frac{
  |\pmi{\varphi}(\zetaprime)|^2}
  {2\sqrt{Q^{(i)}}}
 d\zeta_1 \dotsb \widehat{d \zeta_i} \dotsb d \zeta_n
\tag \num.3
\\
\Cal F^{-1} T \pmi{\phi}(z)
 &=
\frac{1}{(2 \pi)^n}
  \int_{\Bbb R^{n-1}} 
 e^{-\sqrt{-1}(z_1 \zeta_1 + \dotsb + \widehat{z_i \zeta_i} +
      \dotsb + z_n \zeta_n)}
 \tag \num.4
\\
 & \
 \times e^{\mp\sqrt{-1}z_i \sqrt{Q^{(i)}}}
 \frac{{\pmi{\varphi}}(\zetaprime)}
  {2\sqrt{Q^{(i)}}}
 d\zeta_1 \dotsb \widehat{d \zeta_i} \dotsb d \zeta_n.
\endalignat
$$
{}From (\num.4) we have the following:
\proclaim{Lemma~\num}
We write $\Cal F_{\Bbb R^k}$ for
 the Fourier transform  in $\Bbb R^k$ ($k = n-1, n$).
Then,
$$
\alignat1
   \Cal F_{\Bbb R^n}^{-1} T \pmi{\phi}|_{z_i=0}
 &= \frac{1}{4 \pi}
   \Cal F_{\Bbb R^{n-1}}^{-1}
  \left(\frac{\pmi{\varphi}}{\sqrt{Q^{(i)}}}\right)
 (\zprime), 
\\
   \frac{\partial}{\partial z_i}|_{z_i=0}
   \Cal F_{\Bbb R^n}^{-1} T \pmi{\phi}
 &= \frac{\mp{\sqrt{-1}}}{4 \pi}
   \Cal F_{\Bbb R^{n-1}}^{-1}
  \left(\pmi{\varphi} \right)
  (\zprime).
\endalignat
$$
\endproclaim

\redefine\sc{5}
\sec{}
We recall the Fourier transform of the Riesz potential
$$
  \int_{-\infty}^\infty e^{\sqrt{-1} x \xi} \xi_+^\lambda \, d x
  = \sqrt{-1} e^{\frac{\sqrt{-1} \lambda \pi}2} \Gamma(\lambda + 1)
            (\xi + \sqrt{-1} 0)^{-\lambda - 1}
$$
 for a meromorphic parameter $\lambda$.
Letting $\lambda = 0$, we have
$$
   \Cal F^{-1} (\pmi{Y})
 = \frac{\mp \sqrt{-1}}{2 \pi} (z_i \mp \sqrt{-1} 0)^{-1}
    \delta(\zprime),
$$
 where $\delta(\zprime)$ is the Dirac delta function of $n-1$ variables.
Then we have
$$
\alignat1
  \Cal F^{-1} (\pmi{Y} \cdot T \phi)
  &=
  (\Cal F^{-1} \pmi{Y}) * (\Cal F^{-1} T \phi)
\\
  &= \frac{\mp \sqrt{-1}}{2 \pi} (z_i \mp \sqrt{-1} 0)^{-1}
    \delta(\zprime) * f(z)
\\
  &= \frac{\pm 1}{2 \pi\sqrt{-1}} (z_i \mp \sqrt{-1} 0)^{-1} * f(z).
\endalignat
$$
Here the first two convolutions are for $n$ variables $z_1, \dots, z_n$,
 while the last one for only $z_i \in \Bbb R$.
In view of the definition of $\pmi{\phi}$ and $\pmi{f}$,
 we have proved
\proclaim{Lemma~\num}
If $f = \Cal F^{-1} T \phi$ then
 $\pmi{f}=\Cal F^{-1} T \pmi{\phi}$.
\endproclaim
\redefine\sc{6}
\sec{}
We are now ready to complete the proof of Theorem~\ch.2.
By using Lemma~\ch.4 and Lemma~\ch.5,
 the Plancherel formula for $\Bbb R^{n-1}$,
 and the integration formula (\ch.4.3), respectively, 
 we have
$$
\alignat1
 (\pmi{f}|_{z_i=0}, \frac{\partial \pmi{f}}{\partial z_i}|_{z_i=0})
_{L^2(\Bbb R^{n-1})}
&=
 \frac{\pm \sqrt{-1}}{16 \pi^2} 
  (\Cal F_{\Bbb R^{n-1}}^{-1}(\frac{\pmi{\varphi}}{\sqrt{Q^{(i)}}})),
   \Cal F_{\Bbb R^{n-1}}^{-1} (\pmi{\varphi}))
_{L^2(\Bbb R^{n-1})}
\\
&=
 \frac{\pm \sqrt{-1}}{16 \pi^2 (2 \pi)^{n-1}} 
  (\frac{\pmi{\varphi}}{\sqrt{Q^{(i)}}},\pmi{\varphi})_{L^2(\Bbb R^{n-1})}
\\
&=
 \frac{\pm \sqrt{-1}}{2 (2 \pi)^{n+1}} 
 \|\pmi{\phi}\|^2_{L^2(C)}.
\endalignat
$$
Hence,
$$
   (f,f)_W 
  = \frac{1}{2 (2 \pi)^{n+1}}
  \left( \|\phi_+^{(i)}\|^2_{L^2(C)} + \|\phi_-^{(i)}\|^2_{L^2(C)}\right)
  = \frac{1}{2 (2 \pi)^{n+1}}
\|\phi\|^2_{L^2(C)}.
$$
This finishes the proof of the formula (6.2.5) and hence Theorem 6.2.
\qed

\def\sc{7}
\sec{}
The main content of Theorem~\ch.2 is to give
yet another realization of the inner product and of the Hilbert space.
This is very close to the form most known in the case of the
wave equation, where one integrates Cauchy data on the zero time
hypersurface to get the inner product. Note the connection to the
theory of conserved quantities for the wave equation - we end the paper by
making more explicit this final remark: 

When $p = 2$ let us introduce time and space coordinates $(t,x)$ by
$$(t,x) = (t, x_1, \dots , x_{n-1}) = (z_1, \dots , z_n)$$
and the dual variable $k$ to $x$ so that positive-energy solutions
to the wave equation are given by the Fourier transform ($i = \sqrt{-1}$) 
$$u^+(t,x) = \int_{\Bbb R^{n-1}} e^{i(kx - t|k|)} \frac{\varphi^+(k)}{|k|} dk$$
and similarly for negative-energy solutions
$$u^-(t,x) = \int_{\Bbb R^{n-1}} e^{i(kx + t|k|)} \frac{\varphi^-(k)}{|k|} dk$$
where $kx$ denotes the usual scalar product, $|k|$ the Euclidean length,
and $dk = dk_1 \dotsb dk_{n-1}$.
Any solution is the sum of two such: $u = u^+ + u^-$. The energy of the wave
$u$ is given by
$$\Cal E (u) = \frac{1}{2}\int_{\Bbb R^{n-1}} (|u_t|^2 + |\nabla u|^2) dx$$
and is a conserved quantity, i.e. independent of which constant-time hyperplane
we integrate over. It is easy to see that cross-terms drop out, so that on
the Fourier transform side we obtain
$$\Cal E (u) = (2\pi)^{n-1} \int_{\Bbb R^{n-1}} (|\varphi^+(k)|^2 + 
|\varphi^-(k)|^2) dk
$$for the energy. 
Note that this energy only differs from our inner product by a density
factor of $|k|$,
 and that it may be thought of as an integral over $C$. 
In terms of our conformally invariant inner product (6.2.4)
this is up to a constant just 
$$(u, |H|u) = (u^+, Hu^+) - (u^-, Hu^-),$$
where $H = i \partial_t$
is the energy generator (infinitesimal time translations). In the same way, we
have the analogous \lq\lq conserved quantities" for the ultrahyperbolic equation and the
inner product (6.2.4), namely: Let $H_j = i \partial_{z_j}$ be the generator of
translations in the coordinate $z_j$, then for a solution $f$ in the Hilbert space 
$$\Cal E_j(f) = (f, |H_j|f) $$
is invariant under translations in the coordinate $z_j$. Furthermore, the
quantity $\Cal E_j(f)$ can be expressed in terms of an integral of local quantities.
In particular we may use this to study uniqueness and decay properties of
solutions to $\fLap f = 0$. Since this is outside the scope of the present
paper, we shall not do so here; but note that one easy consequence is the fact, 
that if a solution and its normal derivative vanish on a coordinate
hyperplane, then it is identically zero - a classical fact about the energy
(time zero hyperplane).   
       
\redefine\cite{{}}


\widestnumber\key{12}
\Refs
\tolerance = 2000
\ref
   \key \xabste
   \by M\. Abramowitz and I\. A. Stegun (eds.)
   \book Handbook of mathematical functions with formulas, 
   graphs, and mathematical tables
   \finalinfo Reprint of the 1972 edition 
 \publ Dover %
  \yr 1992
%
\endref
\ref
    \key\xbail
    \by W\. N\. Bailey
    \paper Some infinite integrals involving Bessel functions
     \jour  Proc\. London Math\. Soc\.
     \vol 40
      \yr 1935-36
    \pages 37--48
\endref
\ref
    \key \xbz
    \by B\. Binegar and R\. Zierau
    \paper Unitarization of a singular representation of $SO(p,q)$
    \jour Comm\. Math\. Phys\.
    \vol 138
    \yr 1991
    \pages 245--258
\endref
\ref
   \key\Sahi
    \by A. Dvorsky and S. Sahi
    \paper Explicit Hilbert spaces for certain unipotent representations II
    \jour Invent. Math. 
    \vol 138
    \yr 1999
    \pages 203-224
\endref 
\ref
       \key\xerdHigI
       \by A\. Erd\'elyi
       \book Higher Transcendental Functions
       \vol I
       \publ McGraw-Hill
       \publaddr New York
       \yr 1953
\endref
\ref
       \key\xerdIntII
       \by A\. Erd\'elyi
       \book Tables of Integral Transforms
       \vol II
       \publ McGraw-Hill
       \publaddr New York
       \yr 1954
\endref
\ref
        \key\xgs
        \by I\. M\. Gelfand and G\. E\. Shilov     
        \book Generalized Functions, {\rm I}
        \publ Academic Press
        \yr 1964
\endref
\ref
        \key\GSt 
        \by V. Guillemin and S. Sternberg
        \book Variations on a theme by Kepler
        \publ AMS Colloquium Publications, vol. 42 
        \yr 1990
\endref
\ref
        \key\xHo
		\by L\. H\"ormander
		\paper Asgeirsson's mean value theorem and
		related identities
		\jour J\. Func\. Anal\.
		\vol 184
		\yr 2001
		\pages 337--402
\endref
\ref
    \key\xkdecomp
    \by T\. Kobayashi
    \paper  Discrete decomposability of the restriction of
             $A_{\frak q}(\lambda)$
            with respect to reductive subgroups and its applications
    \jour  Invent\. Math\.
    \vol 117
    \yr 1994
    \pages 181--205
; Part II, Ann\. of Math\.  {\bf 147} (1998) 709--729; 
                     Part III, Invent\. Math\. {\bf 131} (1998) 229--256
\endref
\ref
   \key\xkorsI
    \by T\. Kobayashi and B\. \O rsted
    \paper     Analysis on the minimal representation of $O(p,q)$
              --  {\rm I.} Realization via conformal geometry
    \jour preprint
\endref
\ref
   \key\xkorsII
    \bysame %
    \paper     Analysis on the minimal representation of $O(p,q)$
              --  {\rm II.} Branching laws
    \jour preprint
\endref
\ref
     \key \xkos
     \by B\. Kostant
     \paper The vanishing scalar curvature and the minimal unitary
             representation of $SO(4,4)$
     \eds Connes et al
     \inbook Operator Algebras, Unitary Representations, Enveloping Algebras,
                  and Invariant Theory
     \issue      Progress in Math\. %
      \vol 92
      \publ Birkh\"auser
     \yr 1990
      \publaddr Boston
      \pages 85--124
\endref
\ref
   \key \xtodorov
   \by I. T. Todorov 
   \paper Derivation and solution of an infinite-component wave equation
           for the relativistic Coulomb problem
   \inbook  Lecture Notes in Physics,
      Group representations in mathematics and physics
  (Rencontres, Battelle Res. Inst., Seattle, 1969)
   \pages  254--278
   \publ Springer
   \vol 6
    \yr 1970 
\endref
\endRefs
\enddocument